\documentclass[12pt,a4paper,twoside]{amsart}
\usepackage{a4}
\usepackage{amssymb}
\usepackage{graphicx}
\usepackage{amsmath}
\usepackage[english]{babel}
\usepackage[latin1]{inputenc}
\usepackage{lscape}
\usepackage{epstopdf}
\usepackage{enumerate}
\usepackage[left=2.5cm,right=2.5cm,top=3.5cm,bottom=3.5cm]{geometry}
\usepackage[labelfont={footnotesize,sf,bf},textfont={small,sf}]{caption}
\usepackage{multirow}

\newtheorem{defi}{Definition}
\newtheorem{teor}{Theorem}

\begin{document}

\title{Peaks and jumps reconstruction with B-splines scaling functions} 

\author{Luis Ortiz-Gracia}
\address{Centre de Recerca Matem\`{a}tica, 
         Campus de Bellaterra, Edifici C, 
         08193 Bellaterra  
         (Barcelona), Spain}
\email{lortiz@crm.cat}

\author{Josep J. Masdemont}
\address{Departament de Matem\`{a}tica Aplicada I, 
         Universitat Polit\`{e}cnica de Catalunya, 
         Diagonal 647, 
         08028 Barcelona, Spain}
\email{josep@barquins.upc.edu}

\date{October 2012}

\begin{abstract}
We consider a methodology based in B-splines scaling functions to numerically 
invert Fourier or Laplace transforms of functions in the space $L^2(\mathbb{R})$. The original function is approximated
by a finite combination of $j^{th}$ order B-splines basis functions and we provide 
analytical expressions for the recovered coefficients. The methodology is particularly well 
suited when the original function or its derivatives present peaks or jumps due to discontinuities
in the domain. We will show in the numerical experiments the robustness and accuracy of the method. 
\end{abstract}

\maketitle
\renewcommand{\refname}{References}
\renewcommand{\figurename}{Figure}

\section{Introduction}

Fourier and Laplace transforms are useful in a wide number of applications in 
science and engineering. As it is well known, the Fourier transform is
closely related to the Laplace transform for zero value funcions on the
negative time axis.
There is a strong interest in the efficient numerical inversion of Laplace 
transforms (\cite{Abate1996},\cite{Abate2000}) and Fourier transforms, due to 
the fact that the solutions to some problems are known in the transform domain 
rather than in the original domain. An increasing number of papers have recently
appeared to invert Fourier transforms with wavelets, like \cite{Gao2003} and 
\cite{Gao2005} with coiflets wavelets and \cite{Greene2008} with Mexican, 
Morlet, Poisson and Battle-Lemari\'e wavelets. In particular in the Financial 
Engineering context, \cite{Haven2009} inverts a Laplace transform by means of 
B-splines wavelets of order $1$.

Recently, a new method called the COS method developed in 
\cite{Fang2008} for solving the inverse Fourier integral is 
capable to accurately recover a function from its Fourier transform in a 
short CPU time being as well of very easy implementation. However, when the 
function to be recovered presents discontinuities or it is highly peaked, a 
lot of terms in the expansion must be considered to reduce the approximation 
error. 

In this paper we present a novel approximation based on B-splines scaling 
functions to numerically invert Fourier transforms that it is particularly 
suited 
for functions that exhibit peaks and/or jumps in its domain. There are some 
properties that, at first glance, make these basis functions particularly  
suited to approximate such non-smooth functions. B-splines are the 
most regular scaling functions with the shortest support for a given 
polynomial degree. Another important fact is the explicit formulation in the 
time (or space) domain as well as in the frequency domain. However a slight 
drawback of these basis functions is that the system that they form is 
semi-orthogonal, i.e., the scaling functions are orthogonal among different 
scales but not necessarily at the same scale.

In previous work \cite{Masdemont2011} the authors numerically inverted the Laplace transform of a distribution 
function in the interval $[0,1]$ making use of the Haar scaling functions. Now, in the present paper, 
we consider the problem of inverting the Fourier transform to recover an 
$L^2(\mathbb{R})$ function by approximating it by a finite sum of B-splines 
scaling functions of order $j$ (where $j=0$ is the particular case of the Haar 
system). We also provide a list of the different errors accumulating 
within the numerical procedure. So, here, the 
Fourier inversion is carried out with B-splines scaling functions rather than with 
wavelet functions (\cite{Greene2008}) or a combination of both 
(\cite{Gao2003},\cite{Gao2005}). We fix the scale parameter in the wavelet expansion and only remaining is the translation 
parameter, facilitating the inversion procedure. 
 Furthermore, we provide an analytical 
expression for the coefficients of the approximation. 

As will be shown in the section devoted to numerical examples, the Wavelet 
Approximation (WA) that we present is well capable to detect jumps or peaks 
produced by discontinuities in the function itself or in first derivatives.
On the contrary, the COS method is better to approximate analytical functions. 

The paper is organized as follows. In Section 1.1 we give a short literature review regarding the Laplace transform inversion. 
Section 2 gives a brief introduction 
concerning multiresolution analysis and B-splines scaling functions. In Section 3 we present 
the Wavelet Approximation method to recover functions from its Fourier 
transform by means of B-splines scaling functions. Section 4 is devoted to the COS method, while numerical experiments, comparing the 
Wavelet Approximation and the COS method are shown in Section 5. 
Finally, Section 6 concludes.

\subsection{Laplace transform inversion}
Suppose that $f$ is a real- or complex-valued function of the variable $x>0$ and $s$ is a real or complex parameter. We define the \emph{Laplace transform} of $f$ as,
\begin{equation} \label{LT_definition}
\widetilde{f}(s)=\int_0^{+\infty} e^{-sx} f(x)dx=\lim_{\tau \rightarrow +\infty} \int_{0}^{\tau} e^{-sx}f(x)dx,
\end{equation}
whenever the limit exists (as a finite number).

We state the inverse transform as a theorem (see \cite{Dyke2001} for a detailed proof).
\begin{teor} \label{LT_inversion_teor}
(Bromwich inversion integral) If the Laplace transform of $f(x)$ exists, then,
\begin{equation} \label{LT_inversion_formula}
f(x)=\lim_{k \rightarrow +\infty} \frac{1}{2 \pi i} \int_{\sigma-ik}^{\sigma+ik} \widetilde{f}(s)e^{sx}ds, \quad x>0,
\end{equation}
where $\left |f(x) \right | \le e^{\Sigma x}$ for some positive real number $\Sigma$ and $\sigma$ is any other real number such that $\sigma > \Sigma$.
\end{teor}
The usual way of evaluating this integral is via the residues method taking a closed contour, often called the Bromwich contour. 


In this section we present some numerical algorithms to invert the Laplace transform. 
A natural starting point for the numerical inversion of Laplace transforms is the Bromwich inversion integral stated in 
Theorem \ref{LT_inversion_teor}. 
If we choose a specific contour and perform the change of variables $s=\sigma+iu$ in (\ref{LT_inversion_formula}), we obtain an integral of a 
real valued function of a real variable. Then, after algebraic manipulation and after applying the Trapezoidal Rule we obtain,
%
\begin{equation} \label{LT_inversion_formula_real_trap}
f(x) \simeq f_h(x)=\frac{he^{\sigma x}}{\pi} \Re \left(\widetilde{f}(\sigma) \right )+ \frac{2he^{\sigma x}}{\pi} \sum_{k=1}^{+\infty} \Re \left(\widetilde{f}(\sigma+ikh) \right ) \cos(khx),
\end{equation}
where $\Re(z)$ denotes the real part of $z$.
A detailed analysis of the errors can be found in \cite{Abate2000}.

As pointed out in \cite{Abate2000}, the Bromwich inversion integral is not the
only inversion formula and there are quite different numerical inversion algorithms. 
%
We refer the readers to the Laguerre series representation given in \cite{Abate1996}, which is known to be an efficient method for smooth functions. 
However, if $f$ is not smooth at just one value of $x$, then the 
Laguerre method has difficulties at any value of $x$.


In the context of numerical Laplace inversion, \cite{Gao2003} recovers the function $f$ with a procedure based on wavelets.
They consider $s=\beta+i\omega$
in expression (\ref{LT_inversion_formula}), where $\omega$ is a real variable and $\beta$ is a real constant 
that fulfills $f(x)e^{-\beta x} \in L^2 \left ( \mathbb{R} \right)$,
assuming that $f(x)=0$ when $x<0$. Then, equation (\ref{LT_definition}) can be rewritten as,
\begin{equation*}
 \widetilde{f}(\beta+i\omega)=\int_{-\infty}^{+\infty} e^{-\beta x}e^{-i \omega x}f(x)dx.
\end{equation*}
Defining,
\begin{equation} \label{cambio_fh}
 h(x)=f(x)e^{-\beta x}, \; \text{then} \; \hat{h}(\omega)=\widetilde{f}(\beta+i \omega), 
\end{equation}
where $\hat{h}$ denotes the Fourier transform of $h$.
The authors expand $\hat{h}(\omega)$ in terms of Coiflets wavelets,
\begin{equation} \label{hhat_expansion}
 \hat{h}(\omega)=\sum_{k=-\infty}^{+\infty} c_{m,k} \phi_{m,k}(\omega) + \sum_{j=m}^{+\infty} \sum_{k=-\infty}^{+\infty} d_{j,k} \psi_{j,k}(\omega).
\end{equation}
where
$\phi_{m,k}(\omega)=2^{m/2}\phi(2^m \omega -k)$, $\psi_{j,k}(\omega)=2^{j/2}\psi(2^j \omega -k)$, being $\phi$ and $\psi$ the scaling and wavelet functions
respectively.

The next step consists in inverting the expression (\ref{hhat_expansion}) by means of the Fourier inversion formula.
Finally, considering the expression (\ref{cambio_fh}), the formulae of Laplace inversion become,
\begin{equation*}
\begin{split}
 & f_m(x)=\frac{e^{\beta x}}{2^{m+1} \pi} \hat{\phi} \left (-\frac{x}{2^m} \right ) \sum_{k=-\infty}^{+\infty} \widetilde{f} \left (\beta+ i\frac{M_1+k}{2^m} \right ) e^{ixk/2^m}, \\
 & f(x)=\lim_{m \rightarrow +\infty} f_m(x). 
\end{split}
\end{equation*}

One drawback of this approximation is that the wavelet approach involves an infinite product of complex series and the computation of the Fourier transform of 
some scaling functions. This can look intimidating for practical applications and may also take relatively long computational time. 


Based on operational matrices and Haar wavelets, the author in \cite{Chen2001} presents a new method for performing numerical inversion of 
the Laplace transform where only matrix multiplications and ordinary algebraic operations are involved. However, 
the essential step in the method consists in 
expressing the Laplace transform in terms of $\frac{1}{s}$, which is impossible when we just know numerically the transform. 
Another drawback of 
this method is that the matrices become very large for larger scales.

\section{Multiresolution analysis and cardinal B-splines}
A natural and convenient way to introduce wavelets is following the notion of 
multiresolution analysis (MRA). Here we provide the basic definitions and 
properties regarding MRA and B-spline wavelets, for further information 
see \cite{Mallat1989, Chui1992}.
\begin{defi}
A countable set $\{f_n\}$ of a Hilbert space is a Riesz basis if every 
element $f$ of the space can be uniquely written as $f=\sum_n c_nf_n$, and 
there exist positive constants $A$ and $B$  such that,
\begin{equation*}
A\lVert f \rVert^2 \le \sum_n \left|c_n \right|^2 \le B\lVert f \rVert^2.
\end{equation*}
\end{defi}
\begin{defi}
A function $\phi \in L^2(\mathbb{R})$ is called a scaling function, if the 
subspaces $V_m$ of $L^2(\mathbb{R})$, defined by,
\begin{equation*}
V_m=clos_{L^2(\mathbb{R})} \left\{ \phi_{m,k}: k\in \mathbb{Z} \right\}, 
\quad m \in \mathbb{Z},
\end{equation*}
where $\phi_{m,k}=2^{m/2}\phi(2^mx-k)$, satisfy the properties,
\begin{enumerate}[(i)]
\item $\cdots \subset V_{-1} \subset V_0 \subset V_1 \subset \cdots$.
\item $clos_{L^2} \left ( \bigcup_{m \in \mathbb{Z}} V_m \right )=L^2(\mathbb{R})$.
\item $\bigcap_{m \in \mathbb{Z}} V_m=\{0\}$.
\item For each $m$, $\{ \phi_{m,k}: k\in \mathbb{Z} \}$ is a Riesz 
       (or unconditional) basis of $V_m$.
\end{enumerate}
We also say that the scaling function $\phi$ generates a multiresolution 
analysis $\{V_m\}$ of $L^2(\mathbb{R})$.
\end{defi}
 
The $j^\text{th}$ order cardinal $B$-spline function,  $N_j(x)$,
is defined recursively by a convolution:
\begin{equation*}
N_j(x)=\int_{-\infty}^{\infty} N_{j-1}(x-t)N_0(t)dt=\int_0^1 N_{j-1}(x-t)dt, \quad j \ge 1,
\end{equation*}
where,
\begin{equation*}
N_0(x)=\chi_{[0,1)}(x)=\begin{cases} 1 & \text{if } x\in [0,1) \\ 0 & \text{otherwise.} \end{cases}
\end{equation*}
Alternatively,
\begin{equation*}
N_j(x)=\frac{x}{j}N_{j-1}(x)+\frac{j+1-x}{j}N_{j-1}(x-1), \quad j \ge 1.
\end{equation*}
We note that cardinal $B$-spline functions are compactly supported, since the 
support of the $j^\text{th}$ order $B$-spline function $N_j$ is $[0,j+1]$, 
and they have as the Fourier transform,
\begin{equation*}
\widehat{N_j}(w)=\left( \frac{1-e^{-iw}}{iw} \right)^{j+1}.
\end{equation*}

\begin{figure}[ht]
\begin{center}
\begin{tabular}{c}
\includegraphics[scale=0.4,angle=-90]{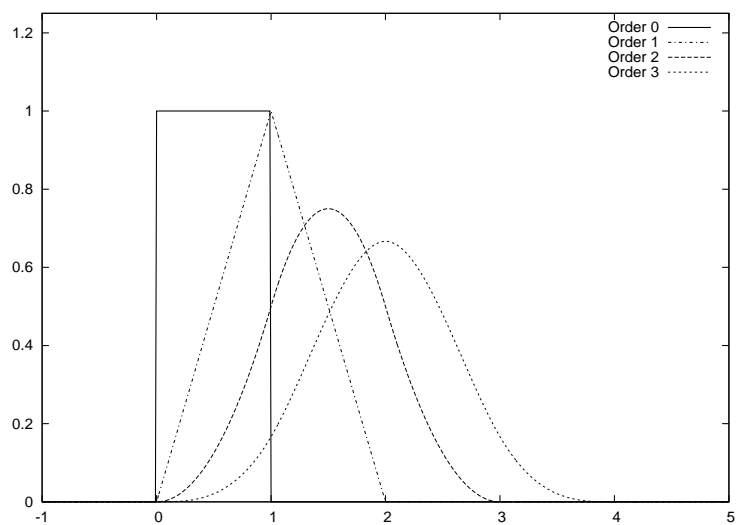}
\end{tabular}
\end{center}
\caption{Cardinal B-splines of orders $j=0,1,2,3$.}
\label{fig-B_splines_0}
\end{figure}

In this paper we consider $\phi^j=N_j$ as the scaling function which generates 
a MRA (see Figure \ref{fig-B_splines_0}). Clearly, for $j=0$ we have the scaling function of the Haar wavelet 
system. We also remark that from the previous discussions, for every function 
$f_m \in V_m$, there exists a unique sequence 
$\{ c_{m,k}^j \}_{k \in \mathbb{Z}} \in l^2(\mathbb{Z})$, such that,
\begin{equation*}
f_m(x)=\sum_{k \in \mathbb{Z}} c_{m,k}^j \phi_{m,k}^j(2^mx-k).
\end{equation*}

\section{The wavelet approximation method}
Let us now consider a function $f \in L^2(\mathbb{R})$ and its Fourier 
transform, whenever it exists:
\begin{equation*}
 \widehat{f}(w)=\int_{-\infty}^{+\infty} e^{-iwx}f(x)dx.
\end{equation*}
Since $f \in L^2(\mathbb{R})$ we can expect that $f$ decays to zero, so it can be well approximated in a finite interval $[a,b]$ by,
\begin{equation}
f^c(x)=
\begin{cases}
f(x) & \text{if $x \in [a,b]$,} \\
0 & \text{otherwise.}
\end{cases}
\end{equation}
Let us approximate $f^c(x) \simeq f^c_{m,j}(x)$ for all $x \in [a,b]$, where,
\begin{equation} 
f^c_{m,j}(x)=\sum_{k=0}^{(j+1)\cdot(2^m-1)} c_{m,k}^j \phi_{m,k}^j \left((j+1) \cdot \frac{x-a}{b-a} \right), \quad j \ge 0,
\end{equation}
with convergence in $L^2$-norm. Note that we are not considering the left and right boundary scaling functions (we refer
the reader to Section 3 in \cite{Maleknejad2010} for a detailed description of scaling functions on a bounded interval).

The main idea behind the Wavelet Approximation method is to approximate 
$\widehat{f}$ by $\widehat{f}^c_{m,j}$ and then to compute the coefficients 
$c_{m,k}^j$ by inverting the Fourier Transform. Proceeding this way,

\begin{equation*}
\begin{split}
 \widehat{f}(w)&=\int_{-\infty}^{+\infty} e^{-iwx}f(x)dx \simeq \int_{-\infty}^{+\infty} e^{-iwx}f^c_{m,j}(x)dx \\
                      & =\sum_{k=0}^{(j+1)\cdot(2^m-1)} c_{m,k}^j \left( \int_{-\infty}^{+\infty} e^{-iwx} \phi_{m,k}^j\left((j+1)\cdot\frac{x-a}{b-a}\right)dx \right). 
\end{split}
\end{equation*}
Introducing the change of variables $y=(j+1)\cdot \frac{x-a}{b-a}$, gives us,
\begin{equation*}
\begin{split}
\widehat{f}(w) &\simeq  \frac{b-a}{j+1} \cdot e^{-iaw} \sum_{k=0}^{(j+1) \cdot (2^m-1)} c_{m,k}^j \int_{-\infty}^{+\infty} e^{-iw\frac{b-a}{j+1}y} \phi_{m,k}^j(y)dy \\
                      &= \frac{b-a}{j+1} \cdot e^{-iaw} \sum_{k=0}^{(j+1) \cdot (2^m-1)} c_{m,k}^j \widehat{\phi}_{m,k}^j \left(\frac{b-a}{j+1} \cdot w \right).
\end{split}
\end{equation*}
Finally, taking into account that $\widehat{\phi}_{m,k}^j(\xi)=2^{-\frac{m}{2}}\widehat{\phi^j}(\frac{\xi}{2^m})e^{-i\frac{k}{2^m}\xi}$ and 
performing the change 
of variables $z=e^{-i \frac{b-a}{2^m (j+1)}w}$, we have,
\begin{equation*}
\widehat{f}\left(\frac{2^m (j+1)}{b-a}i\cdot \log(z) \right) \simeq 2^{-\frac{m}{2}} \frac{b-a}{j+1} \cdot z^{\frac{2^m(j+1)a}{b-a}} \widehat{\phi^j} \left(i\cdot \log(z)\right) \sum_{k=0}^{(j+1) \cdot (2^m-1)} c_{m,k}^j z^k.
\end{equation*}
If we consider, 
\begin{equation*}
P_{m,j}(z)=\sum_{k=0}^{(j+1) \cdot (2^m-1)} c_{m,k}^j z^k \quad \text{and} \quad Q_{m,j}(z)=\frac{2^{\frac{m}{2}} (j+1) z^{-\frac{2^m(j+1)a}{b-a}} \widehat{f}\left(\frac{2^m (j+1)}{b-a}i\cdot \log(z) \right)}{(b-a) \widehat{\phi^j} \left(i\cdot \log(z)\right)},
\end{equation*}
then, according to the previous formula, we have,
\begin{equation} \label{P_approx_Q}
P_{m,j}(z) \simeq Q_{m,j}(z).
\end{equation}
Since $P_{m,j}(z)$ is a polynomial, it is (in particular) analytic inside a disc of the complex plane
$\{z \in \mathbb{C}:\left|z\right|<r\}$ for $r>0$. We can obtain expressions for the
coefficients $c_{m,k}^j$ by means of the Cauchy's integral formula. This is,
\begin{align*}
 c_{m,k}^j &= \frac{1}{2\pi i}\int_{\gamma}\frac{P_{m,j}(z)}{z^{k+1}}dz, \quad k=0,...,(j+1) \cdot (2^{m}-1),
\end{align*}
where $\gamma$ denotes a circle of radius $r$, $r>0$, about the origin.

Considering now the change of variables $z=re^{iu}$, $r>0$, gives us,
\begin{equation} \label{coefs_con_P} 
\begin{split}
c_{m,k}^j & = \frac{1}{2\pi r^k}\int_{0}^{2\pi}\frac{P_{m,j}(re^{iu})}{e^{iku}}du
        = \frac{1}{2\pi r^k}\int_{0}^{2\pi}\left[ \Re(P_{m,j}(re^{iu}))\cos(ku)+\Im(P_{m,j}(re^{iu}))\sin(ku) \right]du, 
\end{split}
\end{equation}
where $k=0,...,(j+1) \cdot (2^{m}-1)$, and $\Re(z)$ and $\Im(z)$ stand for the real and imaginary parts of $z$, 
respectively. 

Note that if $k \ne0$ then we can further expand the expression above by considering,
\begin{equation} \label{coefs_con_P2}
c_{m,k}^j=\frac{2}{\pi r^k}\int_{0}^{\pi}\Re(P_{m,j}(re^{iu}))\cos(ku)du.
\end{equation}

On the other side, since $\widehat{\phi^j} \left(i\cdot \log(z)\right)=\left(\frac{z-1}{\log(z)}\right)^{j+1}$, we have,
\begin{equation}
Q_{m,j}(z)=\frac{2^{\frac{m}{2}} (j+1) z^{-\frac{2^m(j+1)a}{b-a}} \widehat{f}\left(\frac{2^m (j+1)}{b-a}i\cdot \log(z) \right) (\log(z))^{j+1}}{(b-a)(z-1)^{j+1}},
\end{equation}
and it has a pole at $z=1$.
Finally, making use of (\ref{P_approx_Q}) and taking into account the former 
observation, we can exchange $P_{m,j}$ by $Q_{m,j}$ in (\ref{coefs_con_P}) and (\ref{coefs_con_P2}) to 
obtain, respectively,
\begin{equation} \label{trapecios00}
c_{m,0}^j \simeq \frac{1}{2\pi}\int_{0}^{2\pi}\Re(Q_{m,j}(re^{iu}))du,
\end{equation}
and,
\begin{align} \label{trapecios0}
c_{m,k}^j &\simeq \frac{2}{\pi r^k}\int_{0}^{\pi}\Re(Q_{m,j}(re^{iu}))\cos(ku)du, \quad k=1,...,(j+1) \cdot (2^{m}-1),
\end{align}
where $r \neq 1$ is a positive real number. 

In practice, both integrals in (\ref{trapecios00}) and (\ref{trapecios0}) are computed by means of the Trapezoidal Rule, and we can define,
\begin{align*}
 I(k) &=\int_{0}^{\pi} \Re(Q_{m,j}(re^{iu}))\cos(ku)du, \\
 I(k;h) &=\frac{h}{2} \left( Q_{m,j}(r)+(-1)^kQ_{m,j}(-r)+2 \sum_{j=1}^{M-1} \Re(Q_{m,j}(re^{ih_s}))\cos(kh_s) \right),
\end{align*}
where $h=\frac{\pi}{M}$ and $h_s=sh$ for all $s=0,\dots,M$. Proceeding this
way we find,
\begin{equation} \label{coeffs_calculation1}
\begin{split}
 c_{m,k}^j &\simeq \frac{2}{\pi r^k} I(k) \simeq \frac{2}{\pi r^k} I(k;h) \\
         &= \frac{1}{M r^k}\left( Q_{m,j}(r)+(-1)^kQ_{m,j}(-r)+2 \sum_{s=1}^{M-1} \Re(Q_{m,j}(re^{ih_s}))\cos(kh_s) \right), 
\end{split}
\end{equation}
where $k=1,...,(j+1) \cdot (2^{m}-1)$.

Let us summarize four sources of error in our procedure to compute the 
numerical Fourier transform inversion using cardinal B-splines wavelets. These
are: 

\begin{enumerate}[(A)]
\item Truncation of the integration range,
\begin{equation*}
\mathcal{E}_1(x):=f(x)-f^c(x), \quad x \in \mathbb{R} \backslash [a,b].
\end{equation*}
\item The approximation error at scale $m$, 
\begin{equation*}
\mathcal{E}_2(x):=f^c(x)-f^c_{m,j}(x), \quad x \in [a,b].
\end{equation*}
\item The discretization error, which results when approximating the integral $I(k)$ by $I(k;h)$ using the trapezoidal rule. 
We can apply the formula for the error of the compound trapezoidal rule considering,
$$
 q_{m,k}^j(u)=\Re(Q_{m,j}(re^{iu}))\cos(ku), \quad
 \mathcal{E}_3:=I(k)-I(k;h),
$$
and assuming that $q^j_{m,k} \in C^2([0,\pi])$. Then,
\begin{equation} \label{error_typeC}
 \left |\mathcal{E}_3 \right |= \frac{\pi^3}{12M^2} \left| (q_{m,k}^j(\mu))^{\prime \prime} \right|, \quad \mu \in (0,\pi).
\end{equation}
\item The roundoff error. If we can calculate the sum in expression (\ref{coeffs_calculation1}) with a precision of $10^{-\eta}$, 
then the roundoff error after multiplying by a factor $\frac{1}{M r^k}$ is approximately $\mathcal{E}_4:=\frac{1}{M r^k} \cdot 10^{-\eta}$. 
Then, the roundoff error increases when $r$ approaches to 0.
\end{enumerate}

\subsection{Choice of the parameter $r$}
As mentioned before, the choice of the parameter $r$ may have influence in both, the discretization and the roundoff errors. 
In this section
we present a detailed analysis of the errors listed before in order to determine the optimum $r$ value. 
To do this, let us consider $f_1(x)=\chi_{[1/2,1)}(x)$ a step function defined in $[0,1]$, where $\chi$ represents the indicator 
function, and its Fourier transform $\widehat{f_1}(w)=\frac{e^{-iw/2}-e^{-iw}}{iw}$. 

Due to the shape of $f_1$ it seems that 
the best B-splines basis to perform the approximation is based in the Haar scaling 
functions (B-splines of order $0$). In this particular case 
$\Re(Q_{m,j}(re^{iu}))$ and its derivatives up to order $2$ can be computed 
relatively straightforwardly, so that we will be able to calculate the optimal 
value of the parameter $r$ in order to minimize the discretization and the 
roundoff errors. We also demonstrate that the approximation error 
(type (B)) is $0$ in this case. To do so, we consider,

\begin{equation*}
Q_{m,0}(re^{iu})=\frac{r^{2^m} \left( \cos(2^mu)+i\sin(2^mu) \right)-r^{2^{m-1}} \left( \cos(2^{m-1}u)+i\sin(2^{m-1}u) \right)}{2^{m/2}(r\cos(u)+ir\sin(u)-1)},
\end{equation*}
where $u \in [0,\pi]$, and,
\begin{equation*}
\begin{split}
 &\Re(Q_{m,0}(re^{iu}))=\\
&=\frac{r^{2^m+1}\cos(2^mu-u)-r^{2^m}\cos(2^mu)-r^{2^{m-1}+1}\cos(2^{m-1}u-u)+r^{2^{m-1}}\cos(2^{m-1}u)}{2^{m/2}(r^2-2r\cos(u)+1)}.
\end{split}
\end{equation*}
Now, we must choose an appropriate $r$ value to control both the 
discretization and the round-off errors. First of all, we consider the 
discretization error which can be estimated by means of expression 
(\ref{error_typeC}). We note that $q_{m,k}^0 \in C^2([0,\pi])$ 
since,
\begin{equation*}
0<(r-1)^2 \le r^2-2r\cos(u)+1 \le (r+1)^2, \quad \forall u \in [0,\pi], r \neq 1. 
\end{equation*}
So we have,
\begin{equation*}
\begin{split}
 \frac{d}{du}q_{m,k}^0(u)&=\frac{d}{du}\Re(Q_{m,0}(re^{iu}))\cos(ku)-k\Re(Q_{m,0}(re^{iu}))\sin(ku), \\
 \frac{d^2}{du^2} q_{m,k}^0(u)&=\frac{d^2}{du^2}\Re(Q_{m,0}(re^{iu}))\cos(ku)-2k\frac{d}{du}\Re(Q_{m,0}(re^{iu}))\sin(ku)-k^2\Re(Q_{m,0}(re^{iu}))\cos(ku).
\end{split}
\end{equation*}
and,
\begin{equation} \label{ReQ_app}
 \left |\Re(Q_{m,0}(re^{iu}) \right| \le \frac{r^{2^m+1}+r^{2^m}+r^{2^{m-1}+1}+r^{2^{m-1}}}{2^{m/2}(r-1)^2} \simeq r^{2^{m-1}}+o(r^{2^{m-1}}),
\end{equation}
where the last approximation holds for suitably small values of the parameter $r$.
For sake of simplicity, we consider only the terms with smaller exponents in the parameter $r$ for the expressions $\frac{d}{du}\Re(Q_{m,0}(re^{iu}))$ and
$\frac{d^2}{du^2}\Re(Q_{m,0}(re^{iu}))$. Then,
\begin{equation} \label{dReQ_app}
 \left |\frac{d}{du}\Re(Q_{m,0}(re^{iu})) \right | \le \frac{2^{(m-2)/2}r^{2^{m-1}}+A(r)}{(r-1)^4} \simeq r^{2^{m-1}}+o(r^{2^{m-1}}),
\end{equation}
and, 
\begin{equation} \label{ddReQ_app}
 \left |\frac{d^2}{du^2}\Re(Q_{m,0}(re^{iu})) \right | \le \frac{2^{(3m-4)/2}r^{2^{m-1}}+B(r)}{(r-1)^8} \simeq r^{2^{m-1}}+o(r^{2^{m-1}}),
\end{equation}
where $A(r)$ and $B(r)$ are polynomials in $r$ with degree greater than $2^{m-1}$, and the approximations in (\ref{dReQ_app}) and (\ref{ddReQ_app}) 
hold for suitably small values of the parameter $r$. Finally, taking into account 
expressions (\ref{error_typeC}),(\ref{ReQ_app}),(\ref{dReQ_app}) and (\ref{ddReQ_app}) we have,
\begin{equation*}
  \left |\mathcal{E}_3 \right | \le \frac{\pi^3}{12M^2} \left (r^{2^{m-1}}+2kr^{2^{m-1}}+k^2r^{2^{m-1}} \right)+o(r^{2^{m-1}})
=\frac{\pi^3}{12M^2} \left (k^2+2k+1 \right)r^{2^{m-1}}+o(r^{2^{m-1}}).
\end{equation*}
We note that $ \left |\mathcal{E}_3 \right | \rightarrow 0$ as $r \searrow 0$ while the roundoff error is increasing in $r$ as $r$ 
approaches to zero. 

Then, 
the total error should be approximately minimized when the two estimates are equal. This leads to the equation,
\begin{equation*}
 \frac{1}{M r^k} \cdot 10^{-\eta}=\frac{\pi^3}{12M^2} \left (k^2+2k+1 \right)r^{2^{m-1}}.
\end{equation*}
After algebraic manipulation, we find,
\begin{equation} \label{optimal_r}
 r_{m,k}=\left( \frac{12M \cdot 10^{-\eta}}{\pi^3(k^2+2k+1)} \right)^\frac{1}{2^{m-1}+k}, \quad k=0,\dots,2^m-1.
\end{equation}

As mentioned before, the roundoff error arises when multiplying the sum in expression (\ref{coeffs_calculation1}) by the 
pre-factor  $(Mr^k)^{-1}$. 
Let us consider $M=2^m$, then the $k$ of interest is $k=2^m-1$ which is the greatest value that this parameter can take 
(small values of $k$ do not cause roundoff errors). 

\begin{figure}[ht]
\begin{center}
\begin{tabular}{cc}
\includegraphics[scale=0.3,angle=-90]{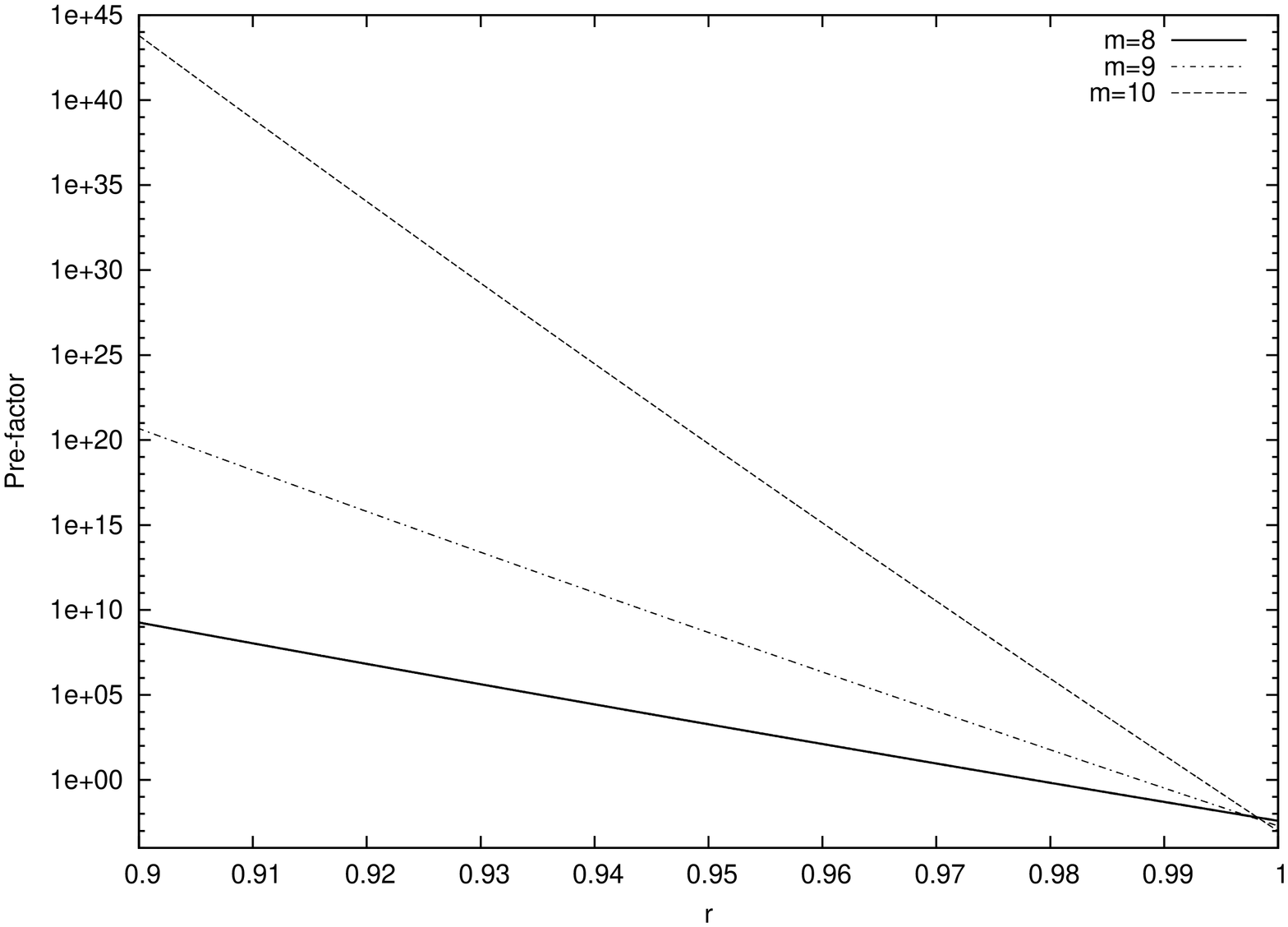} \hfill
\includegraphics[scale=0.3,angle=-90]{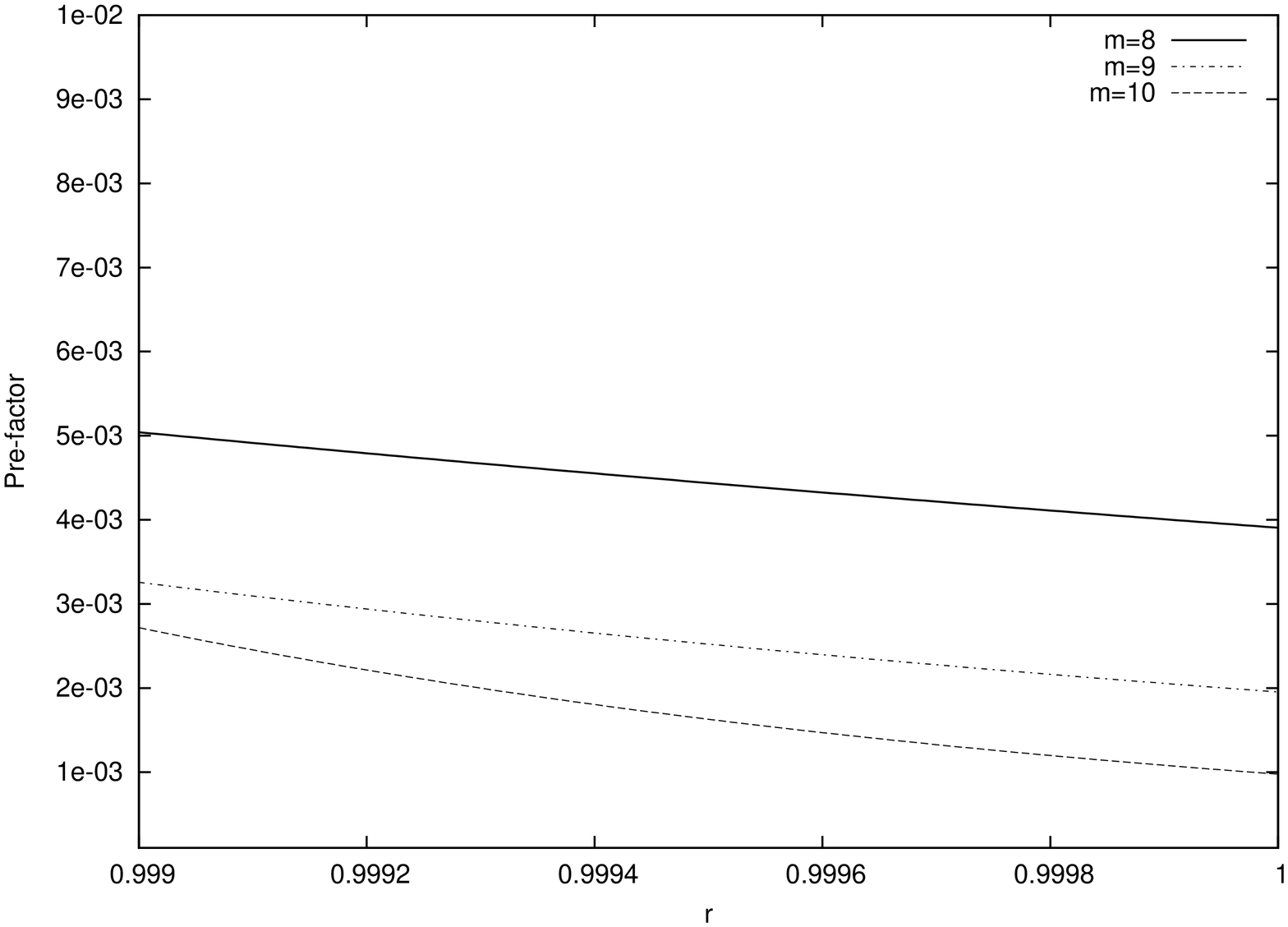}
\end{tabular}
\end{center}
\caption{Pre-factor $(Mr^k)^{-1}$ for $M=2^m$, $k=2^m-1$ and scales $m=8$, $m=9$ and $m=10$.}
\label{fig-prefactor}
\end{figure}

The left plot of Figure \ref{fig-prefactor} represents the pre-factor for values of $r \ge 0.9$, while the right plot shows 
the pre-factor values for $r \ge 0.999$. We also display in Table \ref{parameter-r} the pre-factor values 
$(Mr^k)^{-1}$ for different values of $r$ and scales $m=8$, $m=9$ and $m=10$.

\begin{table}[ht]\footnotesize
\begin{center}
\begin{tabular}{|c| c | c | c|}
\hline
 & \multicolumn{3}{|c|}{Scale} \\
\cline{2-4}
$r$ & $m=8$ & $m=9$ & $m=10$ \\
\hline\hline
$0.9000$ & $1.8194 \cdot 10^9$ & $4.7077 \cdot 10^{20}$ & $6.3040 \cdot 10^{43}$ \\
$0.9100$ & $1.0869 \cdot  10^8$ & $1.6618 \cdot 10^{18}$ & $7.7688 \cdot 10^{38}$ \\
$0.9200$ & $6.6968 \cdot  10^6$ & $6.2395 \cdot 10^{15}$ & $1.0833 \cdot 10^{34}$ \\
$0.9300$ & $4.2522 \cdot  10^5$ & $2.4885 \cdot 10^{13}$ & $1.7047 \cdot 10^{29}$ \\ 
$0.9400$ & $2.7807 \cdot  10^4$ & $1.0529 \cdot 10^{11}$ & $3.0193 \cdot 10^{24}$ \\ 
$0.9500$ & $1.8717 \cdot  10^3$ & $4.7203 \cdot 10^{8}$ & $6.0042 \cdot 10^{19}$ \\ 
$0.9600$ & $1.2960 \cdot  10^2$ & $2.2394 \cdot 10^{6}$ & $1.3373 \cdot 10^{15}$ \\  
$0.9700$ & $9.2550 \cdot 10^{0}$ & $1.1230 \cdot 10^{4}$ & $3.3283 \cdot 10^{10}$ \\   
$0.9800$ & $6.7470 \cdot 10^{-1}$ & $5.9457 \cdot 10^{1}$ & $9.2347 \cdot 10^{5}$ \\ 
$0.9900$ & $5.0674 \cdot 10^{-2}$ & $3.3201 \cdot 10^{-1}$ & $2.8503 \cdot 10^{1}$ \\  
$0.9990$ & $5.0415 \cdot 10^{-3}$ & $3.2566 \cdot 10^{-3}$ & $2.7177 \cdot 10^{-3}$ \\   
$0.9991$ & $4.9145 \cdot 10^{-3}$ & $3.0942 \cdot 10^{-3}$ & $2.4532 \cdot 10^{-3}$ \\ 
$0.9993$ & $4.6699 \cdot 10^{-3}$ & $2.7934 \cdot 10^{-3}$ & $1.9990 \cdot 10^{-3}$ \\ 
$0.9995$ & $4.4376 \cdot 10^{-3}$ & $2.5219 \cdot 10^{-3}$ & $1.6289 \cdot 10^{-3}$ \\ 
$0.9997$ & $4.2169 \cdot 10^{-3}$ & $2.2768 \cdot 10^{-3}$ & $1.3274 \cdot 10^{-3}$ \\ 
$0.9999$ & $4.0072 \cdot 10^{-3}$ & $2.0555 \cdot 10^{-3}$ & $1.0818 \cdot 10^{-3}$ \\ 
\hline
\end{tabular}
\end{center}
\caption{Pre-factor $(Mr^k)^{-1}$ for $M=2^m$, $k=2^m-1$ and scales $m=8$, $m=9$ and $m=10$.}\label{parameter-r} \centering
\end{table}
On the one hand, it is clear that we must concentrate in this second interval in order to get a reasonable roundoff error and 
on the other hand, the discretization error grows when $r$ is close to $1$. Later, in the numerical examples section, we will confirm
the theory developed above for the step function.

Since $f_1$ is compactly supported, we have $\mathcal{E}_1=0$. 

Note that in the case that $j=0$,
\begin{equation*}
 f^c(x)=\sum_{k=0}^{2^m-1}c_{m,k}^0\phi_{m,k}^0 \left(\frac{x-a}{b-a}\right)+\sum_{l=m}^{+\infty}\sum_{k=0}^{2^l-1}d_{l,k}^0\psi_{l,k}^0 \left(\frac{x-a}{b-a}\right),
\end{equation*}
where,
$\{\phi_{m,k}^0\}_{k=0,\dots,2^m-1} \cup \{\psi_{l,k}^0\}_{l \ge m,k=0,\cdots,2^l-1}$ is the Haar basis (orthonormal) system in $L^2 \left([0,1]\right)$. Then,
\begin{equation} \label{error_aproximacion_Haar}
\left \| \mathcal{E}_2 \right \|_{L^2 \left( [a,b] \right)}^2=\left \| f^c-f_m^c \right \|_{L^2 \left( [a,b] \right)}^2 = \left \| \sum_{l=m}^{+\infty}\sum_{k=0}^{2^l-1}d_{l,k}^0\psi_{l,k}\left(\frac{x-a}{b-a}\right) \right \|_{L^2 \left( [a,b] \right)}^2  
=(b-a) \cdot \sum_{l=m}^{+\infty}\sum_{k=0}^{2^l-1}|d_{l,k}^0|^2,
\end{equation}
since $\left \| \psi_{l,k}^0\left(\frac{x-a}{b-a}\right) \right \|_{L^2 \left( [a,b] \right)}^2=(b-a)\cdot \left \| \psi_{l,k}^0\right \|_{L^2 \left( [0,1] \right)}^2=b-a$. Then, the approximation error depends on the length of the interval $[a,b]$ and the detail coefficients,
\begin{equation} \label{detail_coeff}
d_{l,k}^0=\int_{\mathbb{R}} f^c(x) \cdot \psi_{l,k}(x)dx.
\end{equation}

Furthermore, since the detail coefficients (\ref{detail_coeff}) are zero, then we have also that $\mathcal{E}_2=0$ and the approximation 
is exact at any scale level $m$.

\section{The COS method} \label{cos_method}
For completeness, we present here the methodology developed in \cite{Fang2008} 
for solving the inverse Fourier integral\footnote{In order to maintain the 
notation used by the authors, here \begin{equation} \label{FT_Oosterlee} 
\xi(w)=\int_{\mathbb{R}} e^{iwx}f(x)dx, \end{equation} represents the 
characteristic function, and hence the Fourier transform of a density function 
$f(x)$.}. 
The main idea is to reconstruct the whole integral from its Fourier-cosine 
series expansion extracting the series coefficients directly from the 
integrand. Fourier-cosine series expansions usually give an optimal 
approximation of functions with a finite support \cite{Boyd1989}. In fact, the 
cosine expansion of $f(x)$ in $x$ equals the Chebyshev series expansion of 
$f(\cos^{-1}(t))$ in $t$.

For a function supported on $[0,\pi]$, the cosine expansion reads,\\
\begin{equation*}
f(\theta)=\frac{A_0}{2}+\sum_{k=1}^{+\infty} A_k \cos(k\theta),
\end{equation*}
with $A_k=\frac{2}{\pi} \int_0^{\pi}f(\theta)\cos(k\theta)d\theta$. For functions supported on any other finite interval $[a,b] \in \mathbb{R}$, the Fourier-cosine series expansion can easily be obtained via a change of variables,
\begin{equation*}
\theta:=\frac{x-a}{b-a}\pi, \quad x=\frac{b-a}{\pi}\theta+a.
\end{equation*}
It then reads,
\begin{equation*}
f(x)=\frac{A_0}{2}+\sum_{k=1}^{+\infty} A_k \cos(k \pi \frac{x-a}{b-a}),
\end{equation*}
with,
\begin{equation} \label{coeff_approximation}
A_k=\frac{2}{b-a}\int_a^b f(x) \cos(k \pi \frac{x-a}{b-a})dx.
\end{equation}

Since any real function has a cosine expansion when it is finitely supported, the derivation starts with a truncation of the infinite integration range in (\ref{FT_Oosterlee}). Due to the conditions for the existence of a Fourier transform, the integrands in (\ref{FT_Oosterlee}) have to decay to zero at $\pm \infty$ and we can truncate the integration range in a proper way without losing accuracy.

Suppose $[a,b] \in \mathbb{R}$ is chosen such that the truncated integral approximates the infinite counterpart very well, i.e.,
\begin{equation} \label{approximation_CF}
\xi_1(w):=\int_a^b e^{iwx}f(x)dx \simeq \int_{\mathbb{R}} e^{iwx}f(x)dx=\xi(w).
\end{equation} 
Here, $\xi_1$ denotes a numerical approximation.

Comparing equation (\ref{approximation_CF}) with the cosine series coefficients of $f(x)$ on $[a,b]$ in (\ref{coeff_approximation}), we find that,
\begin{equation*}
A_k \equiv \frac{2}{b-a} \Re \left( \phi_1 \left( \frac{k \pi}{b-a} \right) e^{-i \frac{ka \pi}{b-a}} \right),
\end{equation*}
where $\Re$ denotes the real part of the argument. It then follows from (\ref{approximation_CF}) that $A_k \simeq F_k$ with,
\begin{equation*}
F_k \equiv \frac{2}{b-a} \Re \left( \phi \left( \frac{k \pi}{b-a} \right) e^{-i \frac{ka \pi}{b-a}} \right).
\end{equation*}
We now replace $A_k$ by $F_k$ in the series expansion of $f(x)$ on $[a,b]$, i.e.,
\begin{equation} \label{series_cos}
f_1(x)=\frac{A_0}{2}+\sum_{k=1}^{+\infty} F_k \cos(k \pi \frac{x-a}{b-a}),
\end{equation}
and truncate the series summation such that,
\begin{equation} \label{series_cos_trunc}
f_2(x)=\frac{A_0}{2}+\sum_{k=1}^{N-1} F_k \cos(k \pi \frac{x-a}{b-a}).
\end{equation}
The resulting error in $f_2(x)$ consists of two parts, a series truncation error from (\ref{series_cos}) 
to (\ref{series_cos_trunc}) and an error originated from the approximation of $A_k$ by $F_k$. An error analysis that takes 
these different approximations into account is presented in \cite{Fang2008}.

The COS method performs well when approximating smooth functions, but many terms are needed in case that the function or
its first derivative has discontinuities along the domain of approximation.

\section{Numerical examples}
The aim here is to show the accuracy of B-splines to invert Fourier 
transforms of extremely peaked or discontinuous functions with finite support. 
We will consider the B-splines scaling functions of order $j=0,1,2$. We denote by WAi-j the Wavelet Approximation method with 
B-splines of order $i$ at scale $j$, which involves $(2^j(i+1)-i)$ coefficients 
and COS-N the COS approximation method with $N$ terms.

To compute the coefficients of the Fourier inversion in expression (\ref{coeffs_calculation1}), we must set the parameters. 
For this purpose we consider $M=(j+1)\cdot2^m$, where $j$ is the order of the B-spline considered and $m$ is the scale parameter. 
Observe that if we take $M=2k$ instead of $M=(j+1)\cdot2^m$ in (\ref{coeffs_calculation1}) this leads to the following expression,
\begin{equation}
\begin{split}
 c_{m,k} & \simeq \frac{1}{2k r^k}\left( Q_{m,j}(r)+(-1)^kQ_{m,j}(-r)+2 \sum_{s=1}^{2k-1} \Re(Q_{m,j}(re^{is\frac{\pi}{2k}}))\cos\left(ks\frac{\pi}{2k}\right) \right) \\
         &= \frac{1}{2k r^k}\left( Q_{m,j}(r)+(-1)^kQ_{m,j}(-r)+2 \sum_{s=1}^{k-1} (-1)^s \Re(Q_{m,j}(re^{is\frac{\pi}{2k}})) \right),  
\end{split}
\end{equation}
for $k=1,...,(j+1)\cdot(2^{m}-1)$, so the computation time reduces for large scale approximations.

\begin{figure}[ht]
\begin{center}
\begin{tabular}{c}
\includegraphics[scale=0.4,angle=-90]{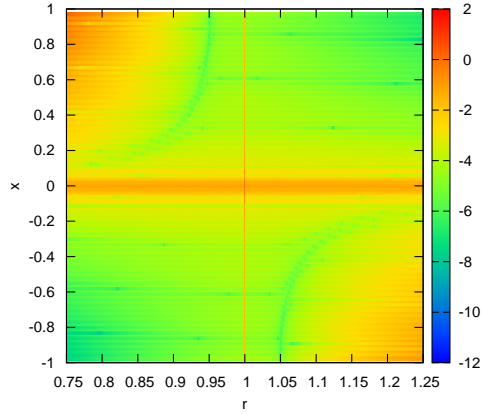}
\end{tabular}
\end{center}
\caption{Logarithm of the absolute error for the approximation of the function $f_2$ ($\alpha=50$) with B-splines of order $1$ at scale $m=5$.}
\label{fig-r_exp_sc5_a50}
\end{figure}

\begin{figure}[ht]
\begin{center}
\begin{tabular}{c}
\includegraphics[scale=0.4,angle=-90]{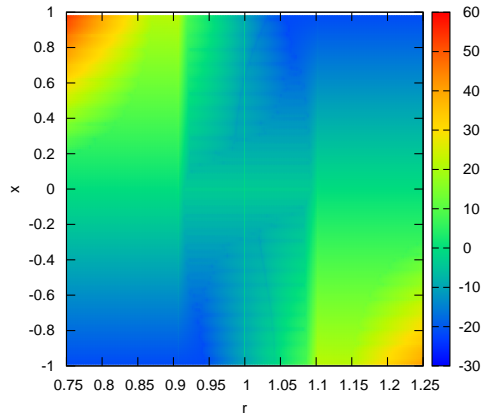}
\end{tabular}
\end{center}
\caption{Logarithm of the absolute error for the approximation of the function $f_2$ ($\alpha=50$) with B-splines of order $1$ at scale $m=9$.}
\label{fig-r_exp_sc9_a50}
\end{figure}

Finally, we set $r=0.9995$. Although we know that $r>0$ and $r \neq 1$, 
we must take into account the two types of errors (C) (discretization) 
and (D) (roundoff) listed previously, which may have 
influence on the computation of the coefficients. Due to the fact that the 
function $\Re(Q_{m,j}(re^{iu}))$ is in general intractable from an analytical 
point of view, we carried out intensive simulations in order 
to understand the influence of parameter $r$.  We 
did these simulations considering the test functions $f_1$ and $f_2(x)=e^{-\alpha |x|}, \alpha >0$ at 
different scale levels and used B-splines scaling functions of order 
$j=0,1,2$. To show just an example, we have plotted the results for $f_2, 
\alpha=50$ and B-splines of order $1$ at scale $m=5$ 
(Figure \ref{fig-r_exp_sc5_a50}) and  $m=9$ (Figure \ref{fig-r_exp_sc9_a50}). 
The colors represent the magnitude of the logarithm of the absolute errors. 
As we can observe, the absolute error remains constant for values 
$|r-1| \le \epsilon$. When $|r-1| > \epsilon$, the error increases 
for high values of $x$ (i.e. high values of the translation parameter $k$) 
and decreases for low values of $x$ (i.e. low values of the translation 
parameter $k$). It is worth mentioning that for high scales, the error grows 
very rapidly when $|r-1| > \epsilon$ (empirical findings demonstrate that 
$\epsilon \simeq 0.05$, for this reason we take approximately the midpoint 
of the interval avoiding $r=1$), while this effect diminishes for shorter 
scales. In fact, as we will see later with the step function, the roundoff error 
almost disappears at very low scales (for instance $m=1$) and the 
discretization error tends to zero when $r$ tends to zero. These facts 
confirm the high impact of the roundoff error at high scale levels and 
only the impact of the discretization error at very low scales.

For comparison, we gave in Section \ref{cos_method} a brief overview of 
the COS method developed by Fang and Oosterlee (\cite{Fang2008}) that is 
based on a Fourier-cosine series expansion and which usually gives optimal 
approximations of functions with finite support (\cite{Boyd1989}). 
The COS method is state of the art and has been applied to efficiently 
recover density functions from their Fourier transforms in order to solve 
important problems arising in Computational Finance.

\begin{figure}[ht]
\begin{center}
\begin{tabular}{cc}
\includegraphics[scale=0.3,angle=-90]{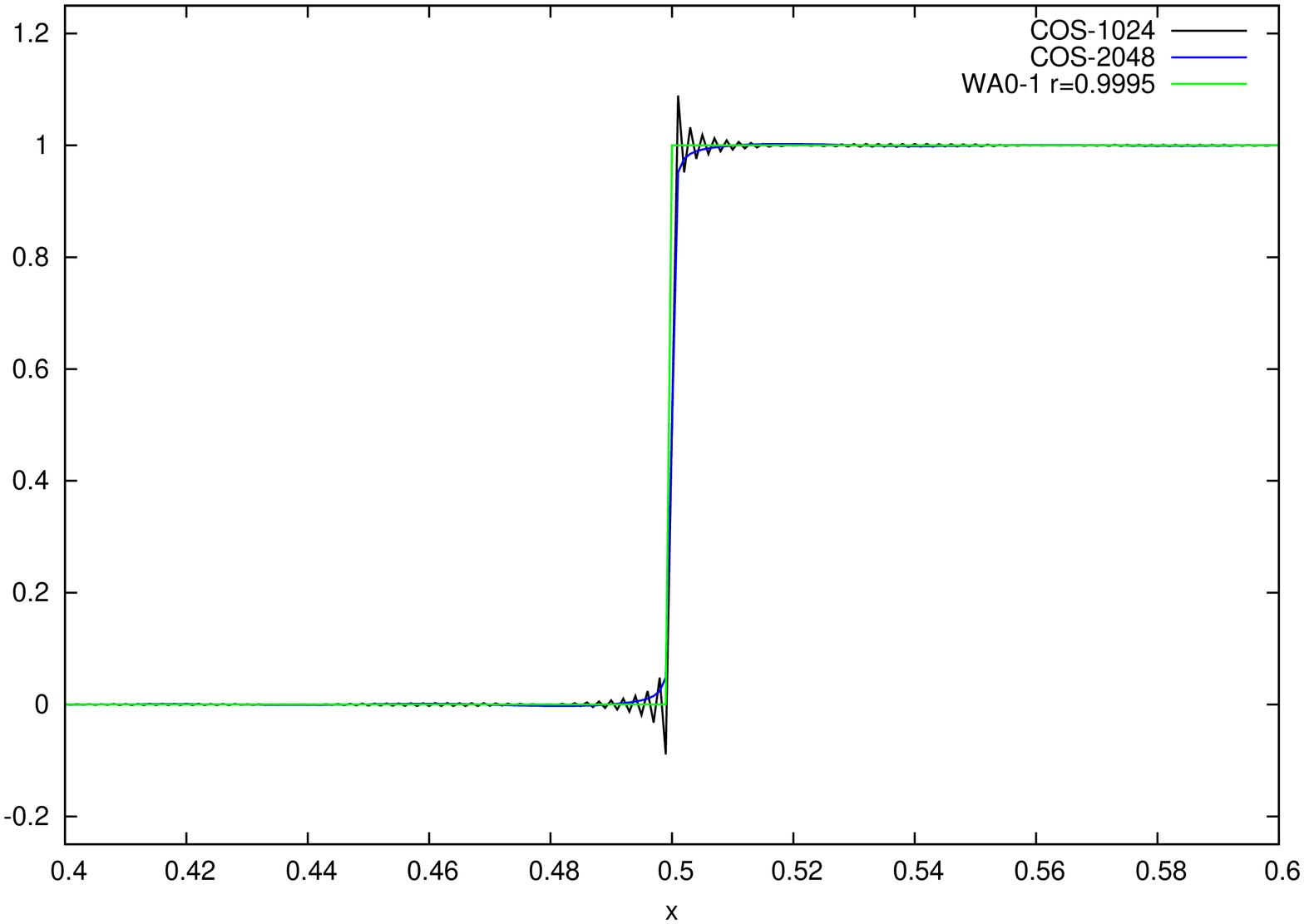}\hfill
\includegraphics[scale=0.3,angle=-90]{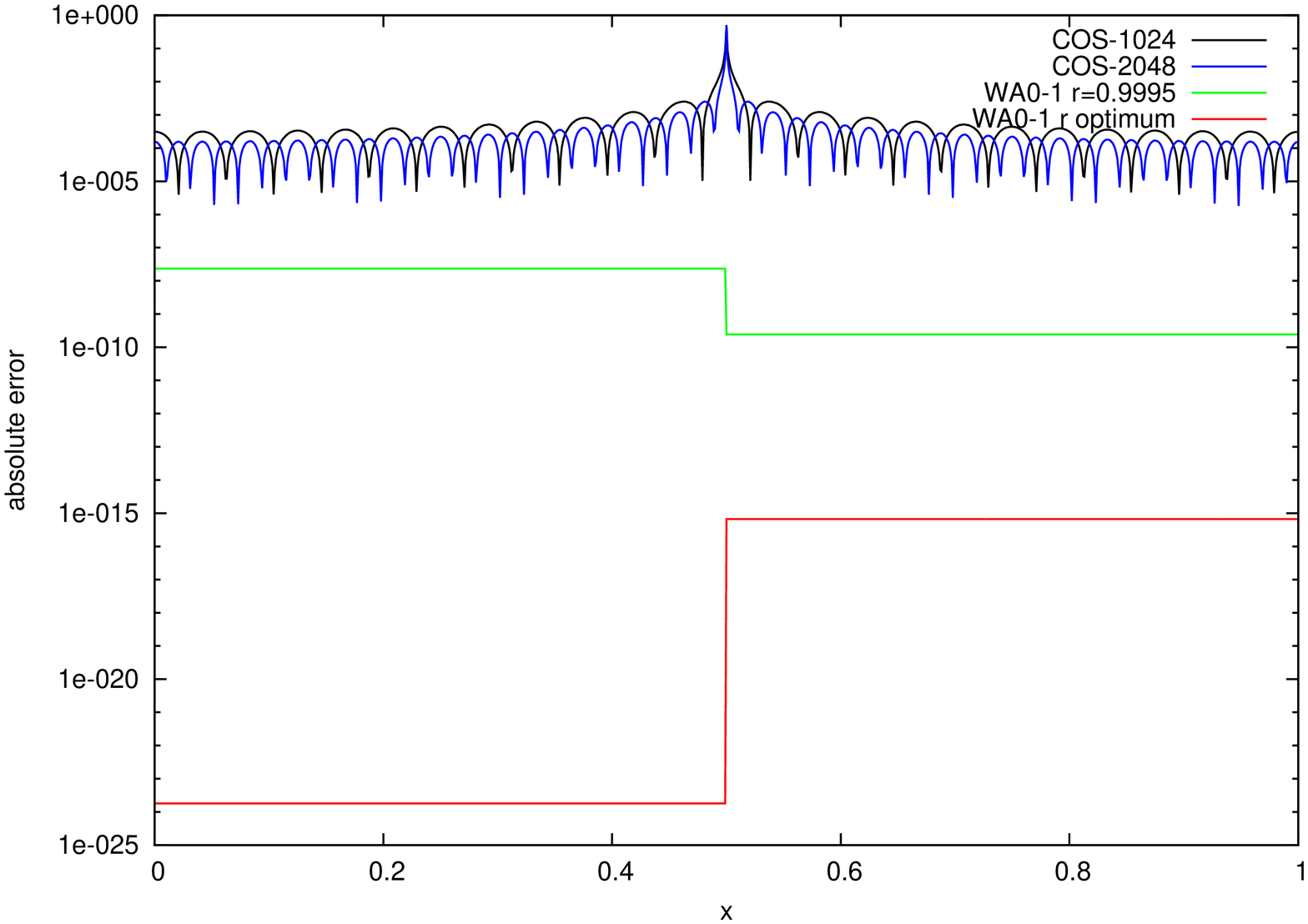} \\ 
\end{tabular}
\end{center}
\caption{Zoom of the approximation (left) and absolute error of the approximation (right) to the function $f_1$.}
\label{fig-f3}
\end{figure}

Let us consider the step function presented before. Observe that if we take $m=1$, there is almost no roundoff error and the only 
remaining error is the discretization error $\mathcal{E}_3$. If we assume that $\eta=16$, then according to (\ref{optimal_r}) 
$r_{1,0}=7.740368 \cdot 10^{-17}$, $r_{1,1}=4.398968 \cdot 10^{-9}$. The left plot of Figure \ref{fig-f3} shows the approximation to 
the step function $f_1$ with the COS and the WA method with $r=0.9995$. The right plot of Figure \ref{fig-f3} represents the 
absolute error of the approximation with the COS and the WA method with $r=0.9995$ and the optimum $r$ computed with expression (\ref{optimal_r}). Note that with just $2$ coefficients the WA method is capable to accurately approximate the step function, while the COS method suffers from the Gibbs phenomenon even if a lot of terms (we show up to $2048$) are added to the COS expansion.

Note that if we consider the approximation to the step function at scale $m=10$, then the values $r_{10,k}$ computed for 
the parameter $r$ are all of them in a neighborhood of $0.95$ in accordance with the massive simulations performed before.

The conclusion in that case is that the COS method performs poorly around the jump and exhibits the Gibbs phenomenon, while Haar
wavelets are naturally capable to deal with this discontinuity. 

\subsection{Exponential function}
\label{ex1}
Let us consider the exponential function, $f_2(x)=e^{-\alpha |x|}, \alpha >0$ and its Fourier transform $\widehat{f_2}(w)=\frac{2\alpha}{\alpha^2+w^2}$.

Observe that this function becomes extremely peaked when $\alpha>>1$. We will test the inversion with $\alpha=50$ and $\alpha=500$. 

\begin{figure}[ht]
\begin{center}
\begin{tabular}{cc}
\includegraphics[scale=0.3,angle=-90]{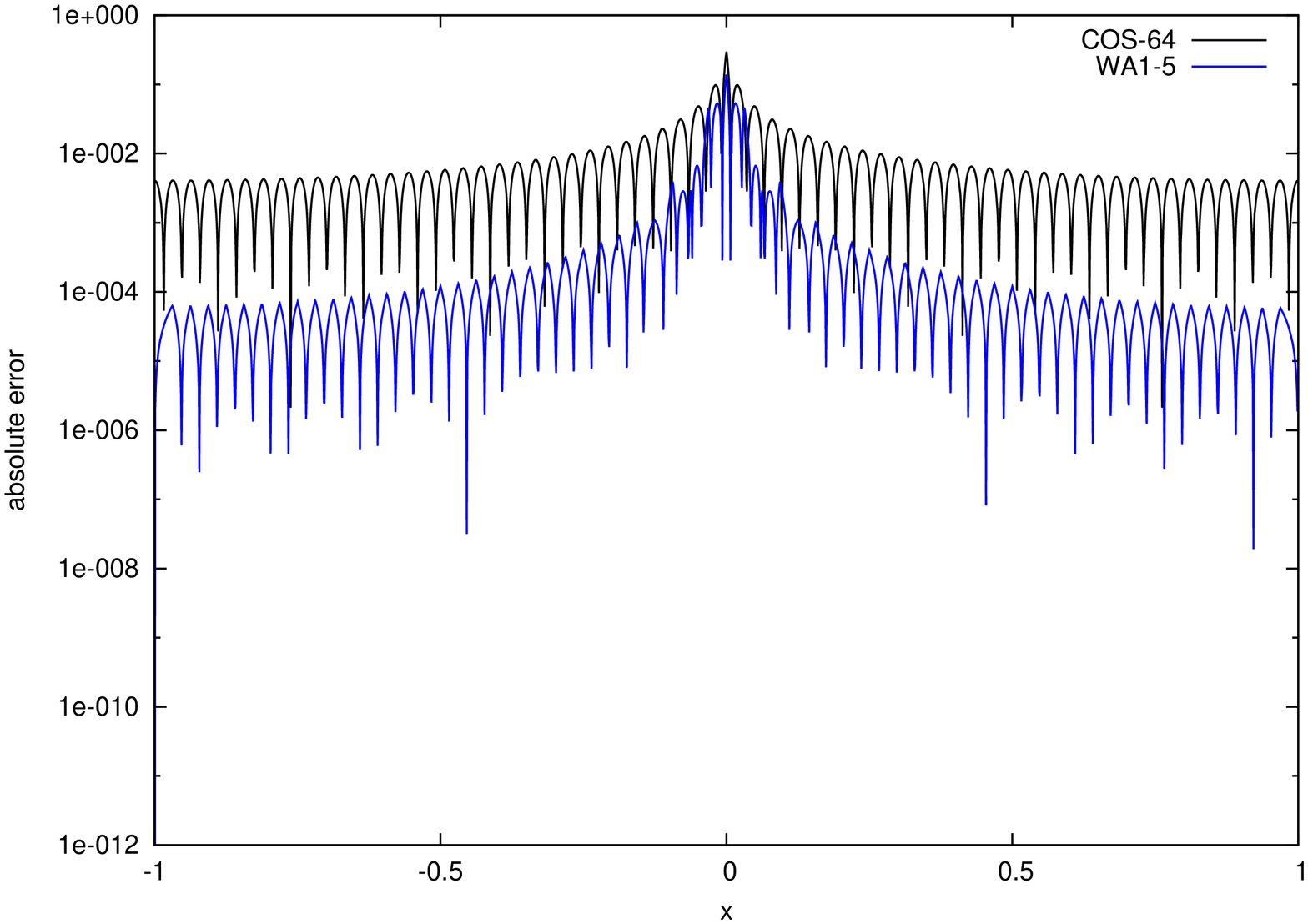}\hfill
\includegraphics[scale=0.3,angle=-90]{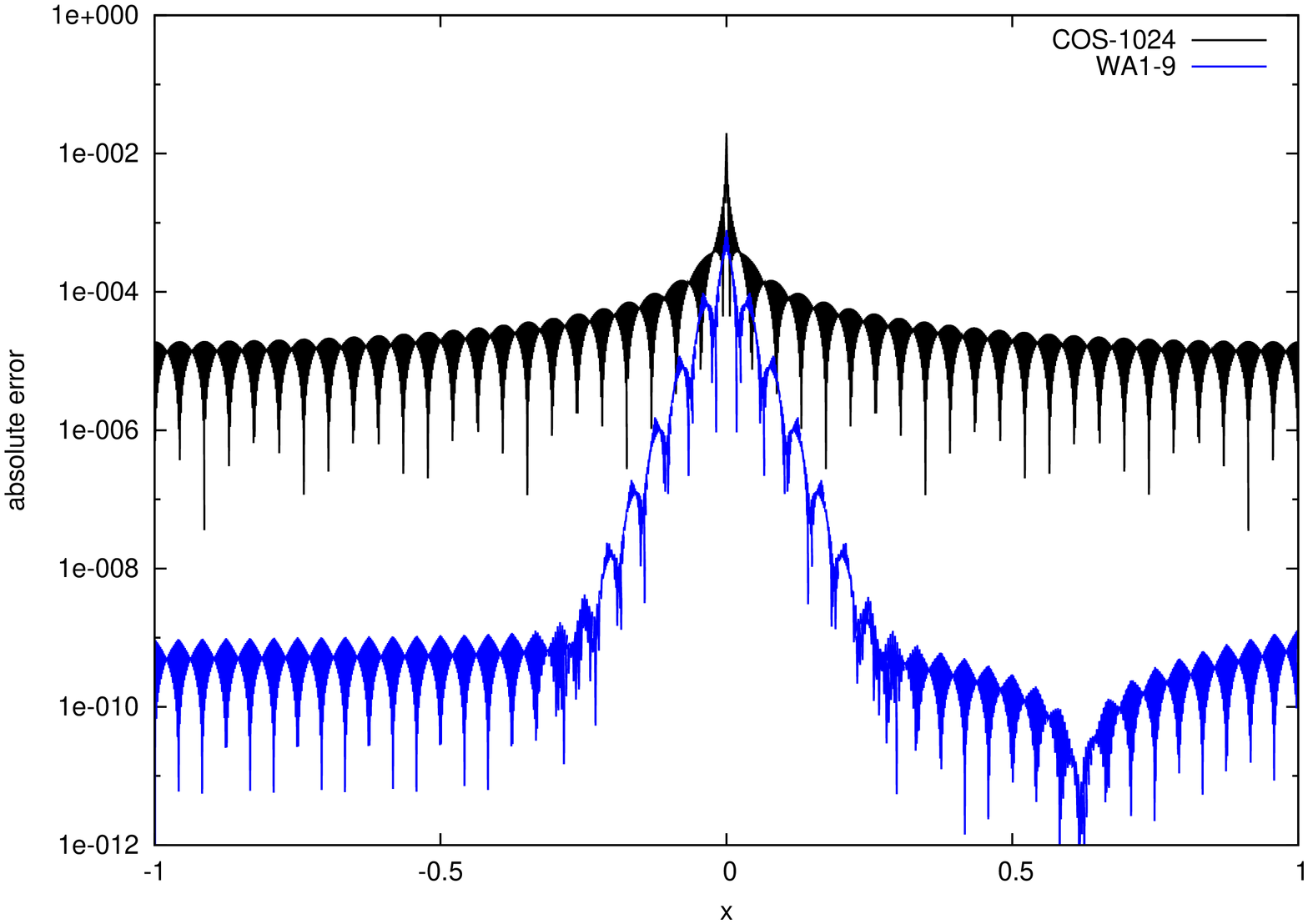} \\ 
\includegraphics[scale=0.3,angle=-90]{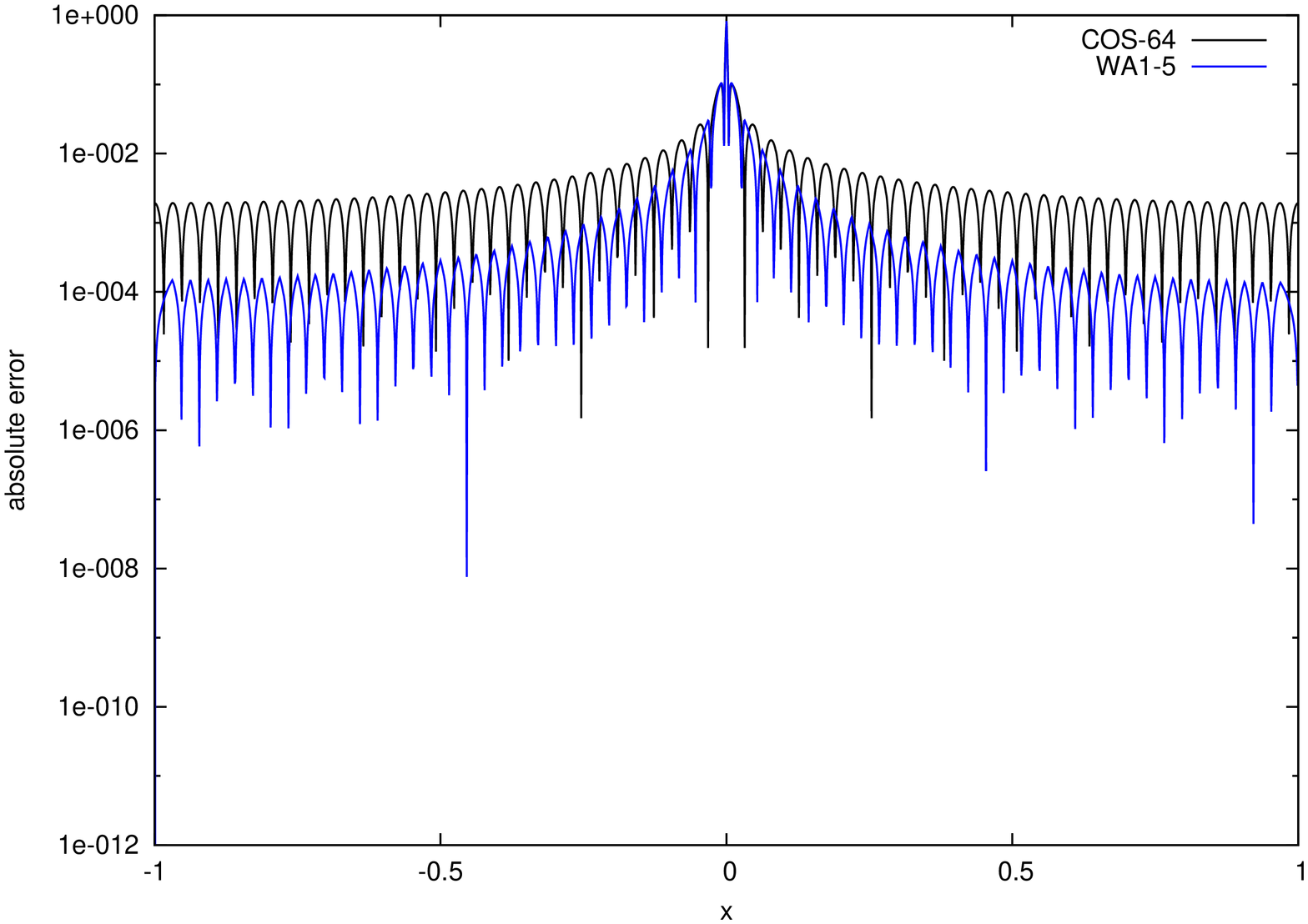}\hfill 
\includegraphics[scale=0.3,angle=-90]{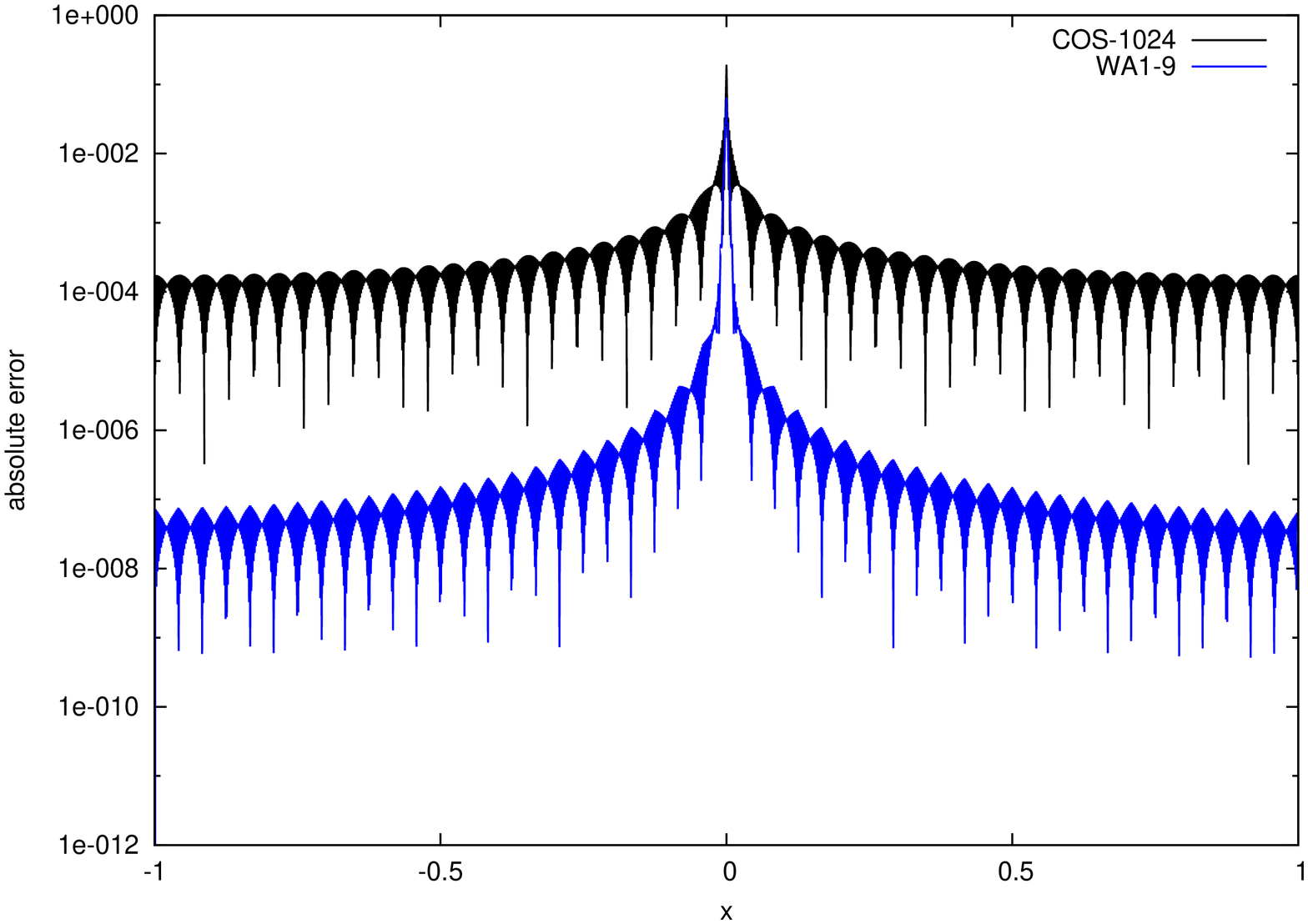} 
\end{tabular}
\end{center}
\caption{Absolute error of the approximation to the function $f_2$ with $\alpha=50$ (top) and $\alpha=500$ (bottom).}
\label{fig-f1}
\end{figure}

We consider the interval $[a,b]=[-1,1]$ 
and $\alpha=50,500$. We carry out the 
comparison between both methods using approximately the same number of terms 
in the approximation. It is worth mentioning that the computation of wavelet 
coefficients is more time consuming than the calculation of COS coefficients. 
Moreover, the COS method is easier to implement. Results are reported in 
Table~\ref{tabla-f1} and plotted in Figure~\ref{fig-f1}. We observe that 
linear B-splines are the most suitable basis functions to approximate highly 
peaked functions like the exponential we are considering. While adding 
many more terms in the COS expansion improves only a bit the approximation, 
when we consider B-splines of order $1$ at higher scales the approximation 
error decays much more faster. The WA method with B-splines of order $1$ 
performs better than the WA method with B-splines of orders $0$ and $2$.

\begin{table}[ht]\footnotesize
\begin{center}
\begin{tabular}{|l| c c| c c|}
\hline  \hline
\multirow{2}{*}{Method}  & \multicolumn{2}{|c|}{$\alpha=50$} & \multicolumn{2}{|c|}{$\alpha=500$} \\ 
& $\min \log|\text{error}|$ & $\max \log|\text{error}|$ & $\min \log|\text{error}|$ & $\max \log|\text{error}|$ \\
\hline \hline
WA0-6 & $-2.928931$ & $-0.374819$ & $-2.523483$ & $-0.040001$ \\ 
WA1-5 & $-7.719564$ & $-0.857977$ & $-8.119854$ & $-0.083319$ \\ 
WA1-9 & $-12.704332$ & $-3.107169$ & $-9.290222$ & $-1.194033$ \\
WA2-4 & $-5.129268$ & $-0.367390$ & $-5.624545$ & $-0.036573$ \\
COS-64 & $-5.673413$ & $-0.526046$ & $-5.825399$ & $-0.057691$ \\
COS-128 & $-5.483920$ & $-0.805937$ & $-5.305269$ & $-0.120147$ \\
COS-256 & $-6.093322$ & $-1.102053$ & $-5.898013$ & $-0.244112$ \\
COS-512 & $-7.488787$ & $-1.402252$ & $-6.433619$ & $-0.450188$ \\
COS-1024 & $-7.454768$ & $-1.703285$ & $-6.500191$ & $-0.716606$ \\
\hline
\end{tabular}
\end{center}
\caption{Approximation errors to the function $f_2$ in the interval $[-1,1]$.}
\label{tabla-f1} \centering
\end{table}

In addition, we consider the following representative examples, 

\subsection{B-spline basis function}
\label{ex3}
We consider the function $f_3(x)=2\phi_{1,0}^1(x)$ in the interval $[a,b]=[0,2]$.

We aim to recover the original function through the Wavelet Approximation inversion method. We apply also the COS method. Like in the previous example, the errors of type (A) and (B) do not have to be considered for this function, since $f_3$ is compactly supported and we carry out the approximation by means of B-splines of order $1$. However, errors of type (C) and (D) remain.

As before, we consider $r=0.9995$. The recovered coefficients are $c_{1,0}^1=2$, $c_{1,1}^1=-5.605010\cdot 10^{-8}$ and $c_{1,2}^1=4.645665\cdot 10^{-8}$, showing high accuracy in the approximation.

\begin{figure}[ht]
\begin{center}
\begin{tabular}{cc}
\includegraphics[scale=0.3,angle=-90]{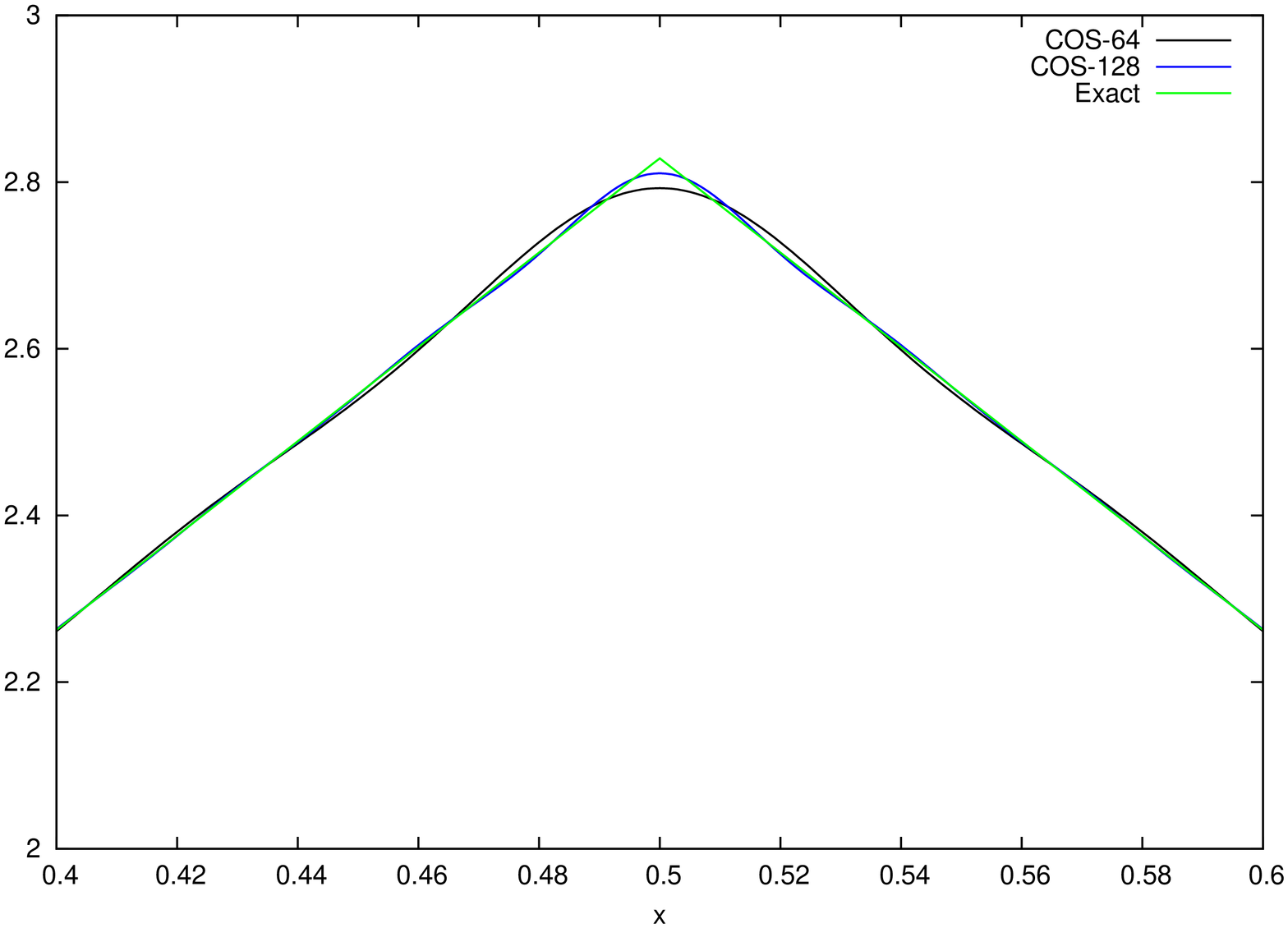}\hfill
\includegraphics[scale=0.3,angle=-90]{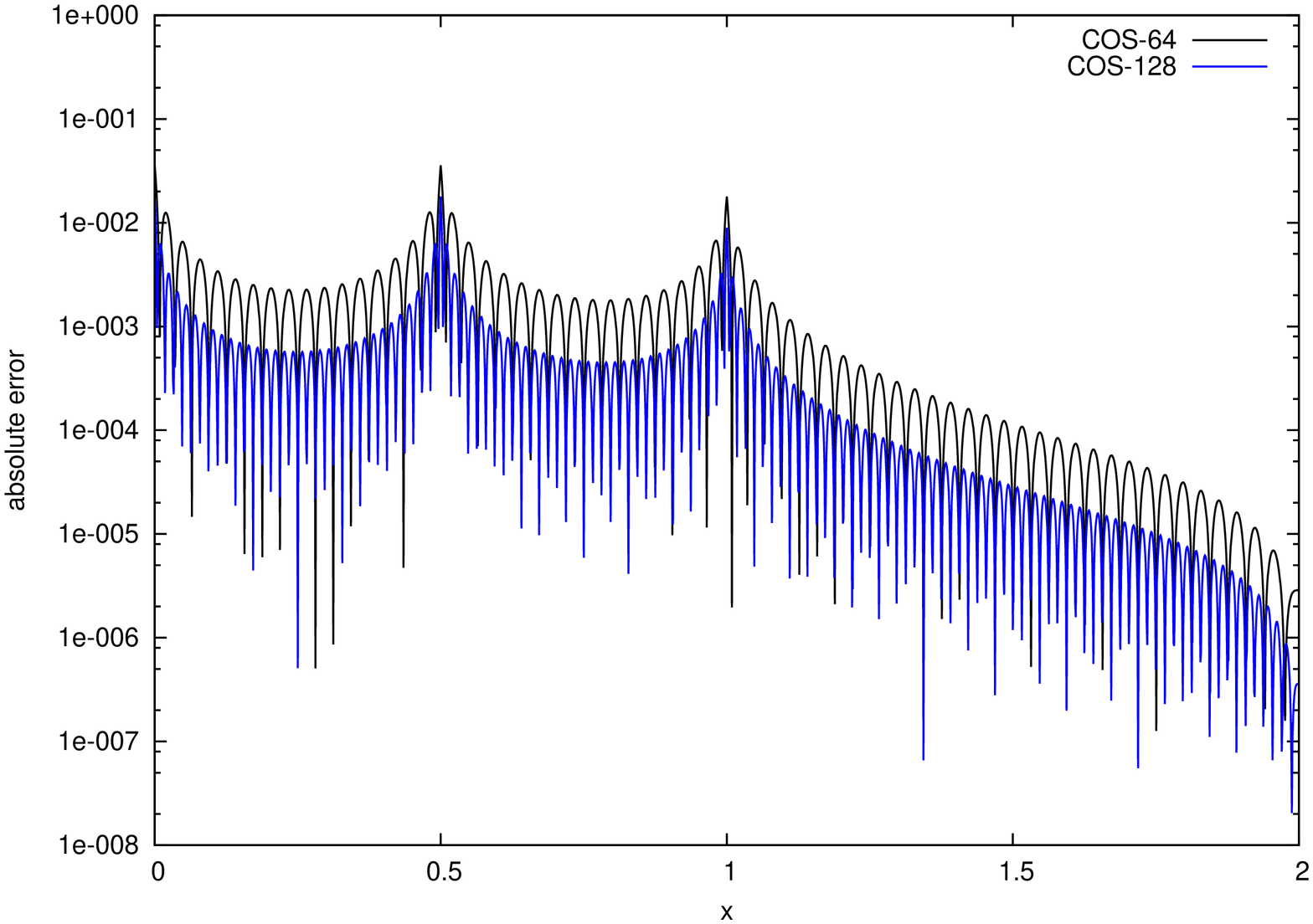} \\ 
\end{tabular}
\end{center}
\caption{Zoom of the approximation (left) and absolute error of the approximation (right) to the function $f_3$.}
\label{fig-f4}
\end{figure}

Figure \ref{fig-f4} shows a zoom of the approximation with COS method to the function $f_3$ in a neighborhood of $0.5$ with $64$ and $128$ terms, while the absolute error is plotted in the right part of the figure. As expected, the error grows up at the non-smooth part of the function, that is, near the points $x=0, x=0.5$ and $x=1$.

\subsection{Combination of B-spline basis functions}
\label{ex4}
We consider the function $f_4(x)=\sum_{k=0}^6 e^{-k} \phi_{2,k}^1(x+1)$ in the interval $[a,b]=[-1,1]$, this is, a finite sum of basis B-splines of order one with coefficients which exhibit exponential decay. Again, the coefficients can be accurately recovered by the Wavelet Approximation method with B-splines of order one at scale $2$. The maximum absolute error is,

$$\max_{k=0\div6}|e^{k}-c_{2,k}^1|=2.681058\cdot10^{-8}.$$

\begin{figure}[ht]
\begin{center}
\begin{tabular}{cc}
\includegraphics[scale=0.3,angle=-90]{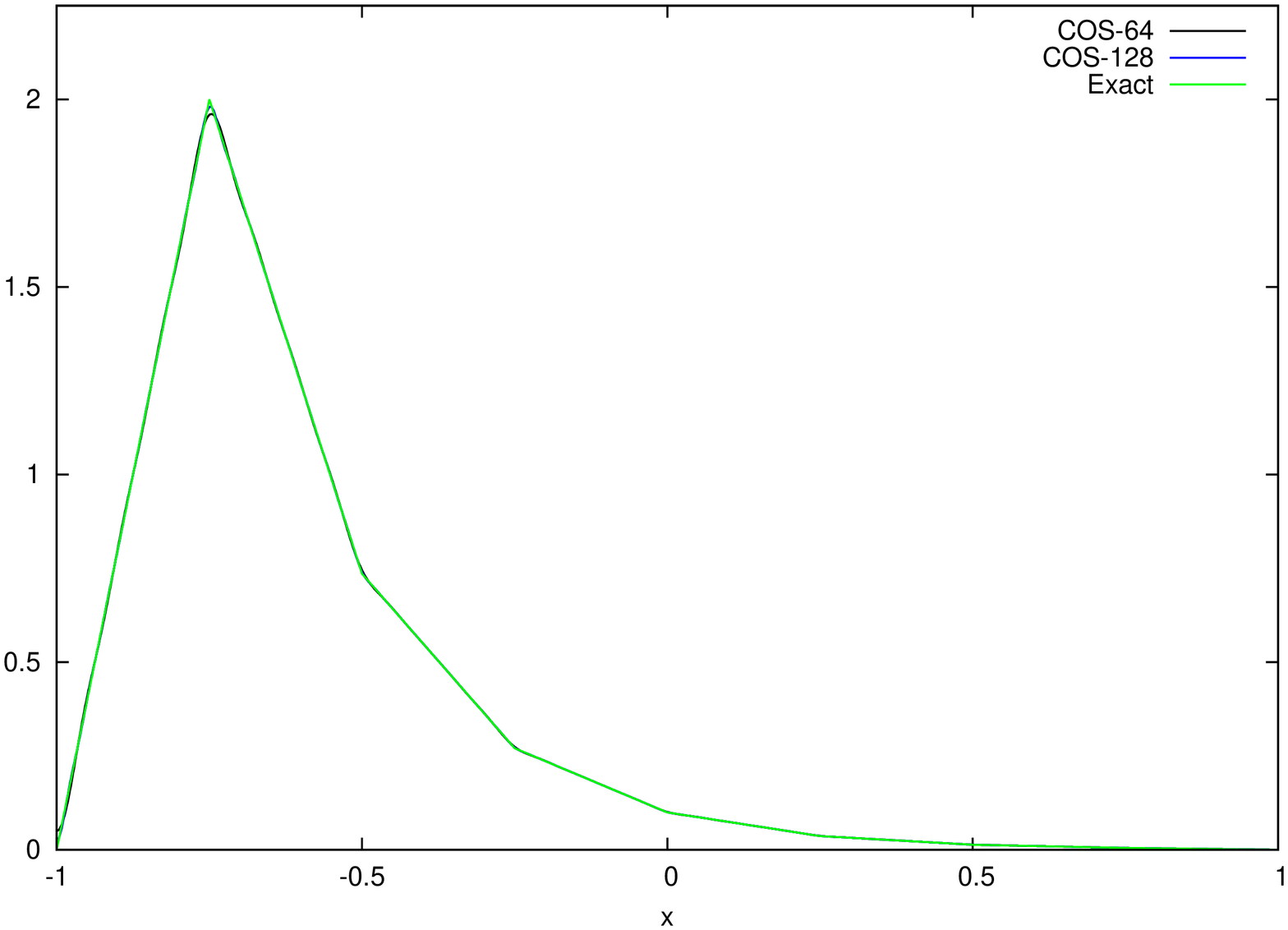}\hfill
\includegraphics[scale=0.3,angle=-90]{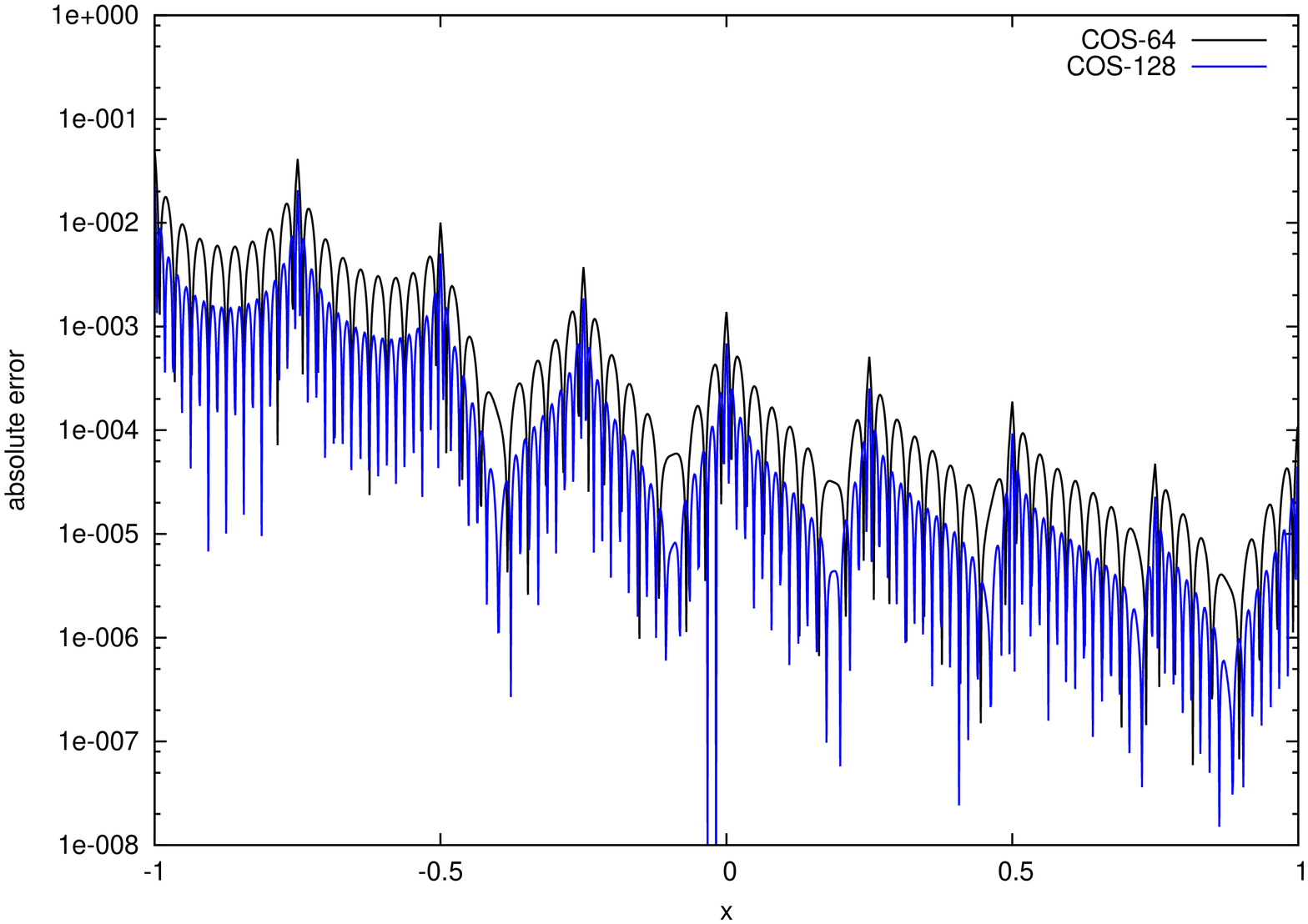} \\ 
\end{tabular}
\end{center}
\caption{Approximation (left) and absolute error of the approximation (right) to the function $f_4$.}
\label{fig-f5}
\end{figure}

Left plot in Figure \ref{fig-f5} represents the approximation to the function with the COS method, while right plot shows the absolute error of the approximation in log scale. As before, the conflictive points for the approximation are the points of non differentiability.

\subsection{Gaussian function}
\label{ex5}
Finally, let us consider the Gaussian function, $f_5(x)=\frac{1}{\sigma \sqrt{2\pi}}e^{-\frac{x^2}{2\sigma^2}}$ and its Fourier 
transform $\widehat{f_5}(w)=e^{-\frac{w^2\sigma^2}{2}}$. 

Since the function $f_5$ is analytic, the COS method performs much more better than the Wavelet Approximation. We can see the results in 
Table \ref{tabla-f5}. As expected, quadratic B-splines behave better than Haar or linear B-Splines, since the maximum error decreases when
the order of the B-Spline increases.

\begin{figure}[ht]
\begin{center}
\begin{tabular}{cc}
\includegraphics[scale=0.3,angle=-90]{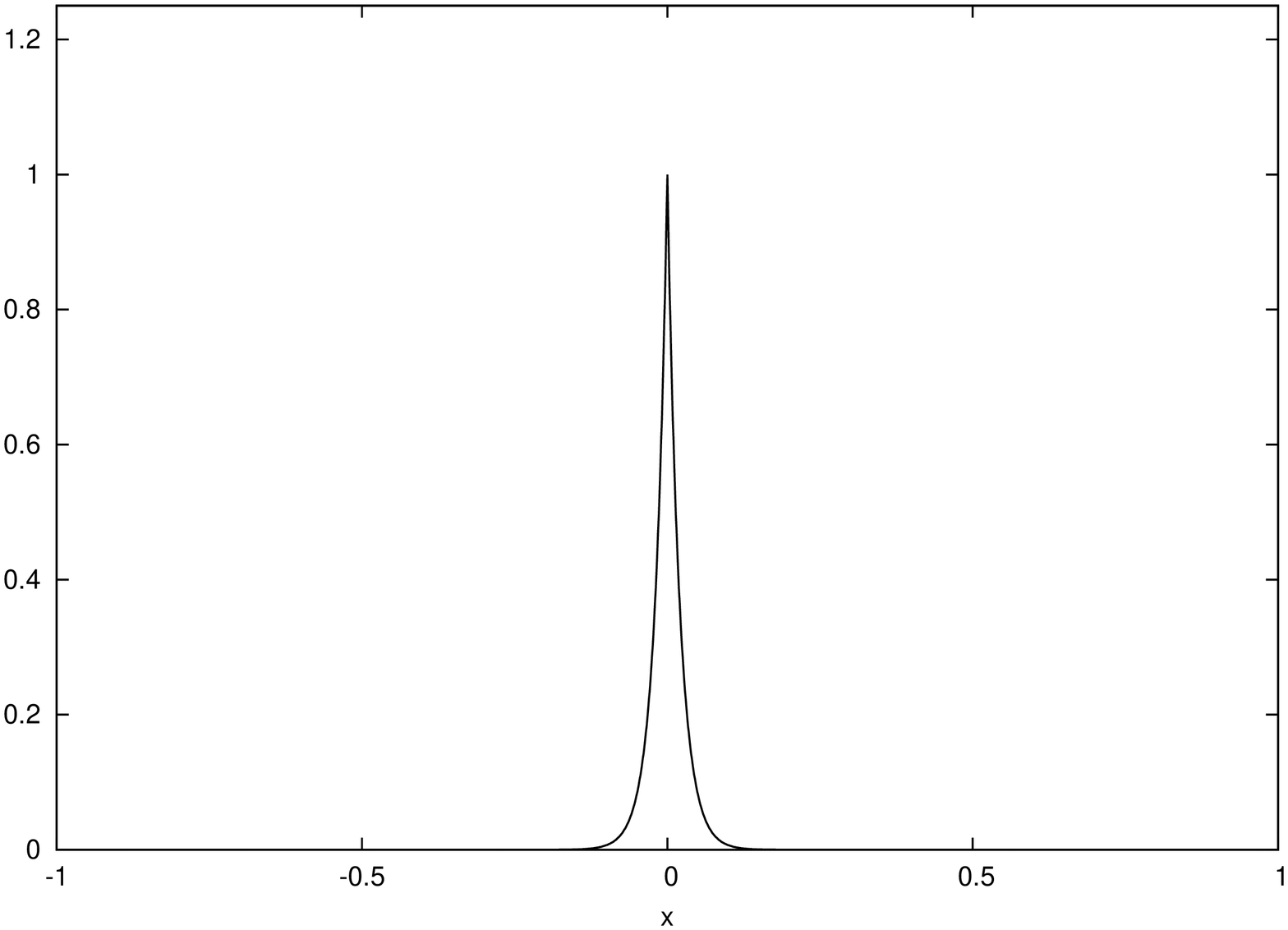}\hfill
\includegraphics[scale=0.3,angle=-90]{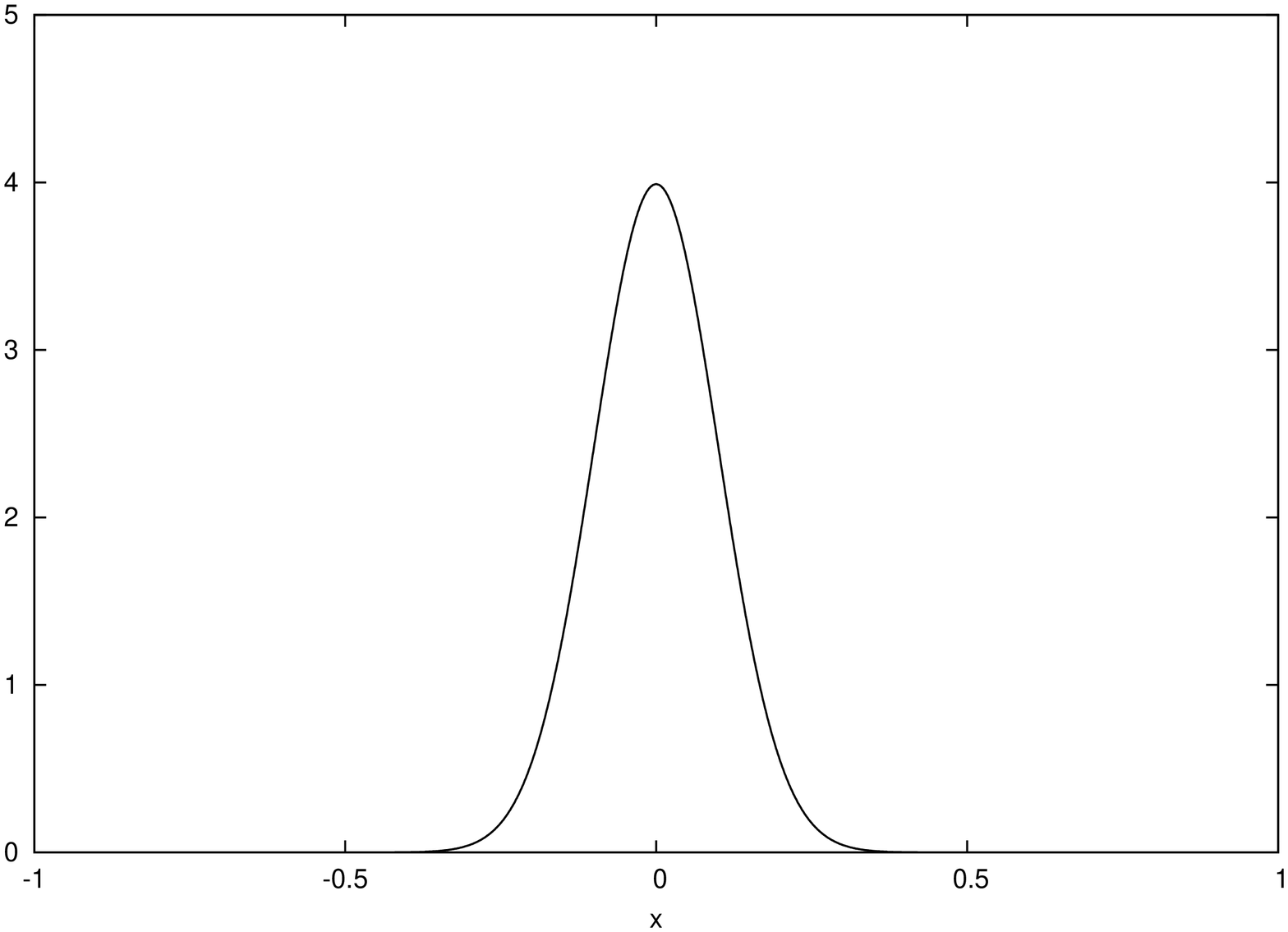} \\ 
\end{tabular}
\end{center}
\caption{The exponential function $f_2$ with $\alpha=50$ (left) and the Gaussian funtion $f_5$ with $\sigma=0.1$ (right).}
\label{fig-f6}
\end{figure}

Left plot in Figure \ref{fig-f6} shows the graph of the exponential function while the right plot represents the graph of the Gaussian
function. It is worth observing the relation between the continuity of the function or its derivatives and the order of the B-Spline
which fits better in the approximation. Haar wavelets for functions with a jump discontinuity, first order B-splines for funtions with 
a jump discontinuity in the first derivative and finally, quadratic B-Splines for more regular functions.

\begin{table}[ht]\footnotesize
\begin{center}
\begin{tabular}{|l| c c|}
\hline  \hline
\multirow{2}{*}{Method}  & \multicolumn{2}{|c|}{$\sigma=0.1$} \\ 
& $\min \log|\text{error}|$ & $\max \log|\text{error}|$ \\
\hline \hline
WA0-6 & $-16.090813$ & $-0.428258$ \\ 
WA1-5 & $-11.687910$ & $-1.482131$  \\ 
WA2-4 & $-9.285967$ & $-2.450365$  \\
COS-32 & $-9.392864$ & $-5.390718$  \\
COS-64 & $-17.699992$ & $-7.530656$ \\
\hline
\end{tabular}
\end{center}
\caption{Approximation errors to the function $f_5$ in the interval $[-1,1]$.}
\label{tabla-f5} \centering
\end{table}

\section{Conclusions}
We have investigated a numerical procedure, the Wavelet Approximation method, to invert the Fourier transform of functions 
with finite support by means of the B-splines scaling functions. First of all, we truncate the function in a sufficiently 
wide interval and then we approximate it by a finite combination of B-splines wavelets up to order $2$. 
Finally the function is recovered from its Fourier transform. 

We have tested and compared the accuracy of this method versus the COS method for a set of heterogeneous functions in terms of the 
continuity and smoothness. 
We have considered a continuous function with a jump discontinuity in its first derivative and a function with a 
jump discontinuity in its domain of definition. Additionally, combinations of B-splines of order one have been also considered as 
the functions to be recovered. With these last two examples we are able to assess the discretization and roundoff errors of the 
Wavelet Approximation method, since no errors of type (A) and (B) take place. Finally, we also have considered the Gaussian function. 
COS method performs better with infinitely differentiable functions, like the Gaussian, due to the fact that the coefficients of 
the expansion decay exponentially fast. On the contrary, WA method is more suitable for functions that exhibit peaks or discontinuities 
along its domain. For the function with a jump discontinuity, B-splines of order $0$ (Haar) are better, while B-splines of order $1$ 
fit accurately the continuous function with a jump discontinuity in its first derivative. Furthermore, little improvement is 
achieved for these two last functions, when adding 
a lot of terms in the COS expansion.   

Due to the fact that COS coefficients are easier to compute in terms of CPU time, hybrid methods involving COS and WA, or 
simply combinations of the WA method at different scales and/or with different orders of the B-splines, could be investigated in the future. 
B-splines are a semi-orthogonal system, so it would be also advisable to explore the orthogonal Daubechies or Battle-Lemari\'e scaling 
functions. Furthermore, it is well known that  standardized B-splines tend to the Gaussian function when the order tends to infinity, 
so they could be very useful to recover Gaussian functions from its Fourier transform.

\section*{Acknowledgements}
This work has been partially supported by the grants SGR2009-859, MTM2009-06973
and MTM2012-31714.

\end{document}